%% file: main.tex
\documentclass[letterpaper,12pt]{article}
\usepackage{braket}
\usepackage{diagbox}
\usepackage{pdflscape} 
\usepackage{adjustbox}

\usepackage{tabularx} 
\usepackage{amsmath,amsthm,amssymb}  
\usepackage{amsmath}
\usepackage{amsfonts}
\usepackage{amssymb}
\usepackage{graphicx} 
\usepackage[margin=0.8in,letterpaper]{geometry} 
\usepackage[english]{babel}
\usepackage{caption}
\usepackage{float}
\usepackage[table]{xcolor}
\usepackage{array}
\usepackage{algorithm}
\usepackage{algpseudocode}
\usepackage{cases}
\usepackage{titlesec}
\usepackage{tikz}
\usepackage[utf8]{inputenc}
\usepackage[T1]{fontenc}
\usepackage{bm}
\usepackage{mathtools}
\usepackage{extarrows}
\usepackage{booktabs} 
\usepackage{dsfont}
\usepackage{authblk}
\usepackage{multirow}
\usepackage{hyperref}
\usepackage{parskip}

\usepackage{rotating}
\usepackage{csquotes}

\newtheorem{theorem}{Theorem}
\newtheorem{lemma}{Lemma}
\newtheorem{proposition}{Proposition}
\newtheorem{corollary}{Corollary}
\newtheorem{definition}{Definition}

\newtheorem{remark}{Remark}

\usepackage{microtype}      
\usepackage{lipsum}
\usepackage{nicefrac}
\usepackage{fancyhdr}       

\author[1]{El-Mehdi Mehiri\\
\texttt{mehiri314@gmail.com}}

\title{\textbf{Explicit M-Polynomial and Degree-Based Topological Indices of Generalized Hanoi Graphs}}

\date{}

\begin{document}

\maketitle

\begin{abstract}
\noindent The M-polynomial, introduced by Deutsch and Klavžar in 2015, provides a unifying algebraic framework for the computation of numerous degree-based topological indices such as the Zagreb, Randić, harmonic, and forgotten indices. Despite its broad applications in chemical graph theory and network analysis, closed expressions of the M-polynomial remain unknown for many important graph families. 

\noindent In this work we derive, for the first time, a complete explicit expression of
the M-polynomial of the generalized Hanoi graphs $H_p^n$ for arbitrary positive 
$p$ and $n$.
Our derivation relies on a detailed combinatorial analysis of the 
occupancy-based structure of $H_p^n$, refined using Stirling and 
$2$-associated Stirling numbers to enumerate all configurations with prescribed 
singleton and multiton counts.
We obtain closed formulas for all diagonal and off-diagonal coefficients of the 
M-polynomial and show how these expressions yield exact values of the main 
degree-based topological indices. 
The correctness of the formulas is supported   through  numerical computation in small instances.
These results provide a complete degree-based description of $H_p^n$ and make
their structural complexity fully accessible through the M-polynomial framework.

\end{abstract}



\section{Introduction}

Topological indices are numerical graph invariants designed to encode structural
information in a compact algebraic form.
Originally motivated by chemical graph theory, where molecular structures are
modeled by graphs whose vertices represent atoms and edges represent bonds,
these indices have found applications in quantitative structure--property and 
structure--activity relationships (QSPR/QSAR); see 
\cite{gutman2012mathematical, Todeschini2010, trinajstic2018chemical}.
Among the numerous classes of such invariants, {\em degree-based} indices occupy
a central position due to their simplicity, interpretability, and strong empirical
correlations with a wide range of physicochemical properties \cite{Dearden2017}.

Such indices have demonstrated strong correlations with a wide range of molecular attributes, including heat of formation, boiling point, strain energy, rigidity, and fracture toughness. As a result, topological indices play a vital role in cheminformatics, drug design, and materials science~\cite{rucker1999topological,klavvzar1996comparison,DENG20113017,bentham,ZHANG1996147,shanmukha2020degree}.

The last decade has seen a growing emphasis on  unifying frameworks  that
encode entire families of topological indices within a single polynomial.
For distance-based indices, this role is played by the Hosoya and Wiener
polynomials \cite{Hosoya1988, Wiener1947}. A wide range of such indices can be efficiently derived from the Hosoya polynomial, as demonstrated in~\cite{Computing}.

In the degree-based setting, Deutsch and Klavžar introduced in 2015 the
\emph{M-polynomial}~\cite{EmericMploy}, defined by
\[
    M(G;x,y)=\sum_{i\le j} m_{i,j}(G)x^iy^j ,
\]
where $m_{i,j}(G)$ counts the number of edges whose endpoints have degrees
$\{i,j\}$. This polynomial serves as a unifying framework for expressing many degree-based topological indices in closed form, much like the Hosoya polynomial does for distance-based indices~\cite{kwun2017m,Chaudhry03102021}. 
Once $M(G;x,y)$ is known, a broad spectrum of classical indices---including the 
Zagreb, Randić, harmonic, symmetric division, forgotten, and inverse sum 
indegree indices---can be obtained by applying differential and integral 
operators to $M(G;x,y)$.

The strength of the M-polynomial lies in its comprehensive encoding of degree-related structural information. It covers data derived from the vertex degrees and the local symmetries of molecular graphs, particularly those representing 2D molecular lattices \cite{Raheem02042020}. As such, the M-polynomial facilitates the computation of a broad spectrum of graph invariants that are often correlated with the physical and chemical properties of the materials under study.

Despite this unifying strength, explicit closed forms of the M-polynomial remain
unknown for many important graph families, particularly those featuring
self-similarity, recursive structure, or constrained dynamics.
A notable example is the family of \emph{generalized Hanoi graphs} $H_p^n$,
which arise from the Tower of Hanoi problem with $p\ge 3$ pegs and $n$ 
distinctly sized discs. Generalized Hanoi graphs occupy an important position among self-similar and
recursively defined graphs; for about these graphs,  
see the monograph of Hinz et al.~\cite[Chapter 5]{hinz2018tower}.
Vertices correspond to legal states of the puzzle, while edges correspond to
legal moves. 
Unlike the Sierpiński graphs $S_p^n$, which form a close sibling family and whose
M-polynomials have been studied in \cite{ishfaq2024topological,Fan,JahangeerBaig03102019}  (See~\cite[Chapter 4]{hinz2018tower} for more about Sierpiński Graphs) ,   the M-polynomial of $H_p^n$ has not been computed
for any $p\ge 4$ so far.

The classical case $p=3$ (the standard Tower of Hanoi graph) has been treated in
\cite{gaur2020m,Asmat}, but the transition from three to four
pegs dramatically increases combinatorial complexity.  This makes the derivation of the M-polynomial for $p\ge 4$ a nontrivial problem.

This paper provides the first explicit closed-form expression of the 
M-polynomial of generalized Hanoi graphs $H_p^n$ for arbitrary 
$p\ge 3$ and $n\ge 0$.
Our main contributions are as follows:
\begin{itemize}
    \item We develop a complete \emph{degree-based combinatorial analysis} of
          the graph $H_p^n$, using the fact that the degree of a vertex depends solely
          on the number of occupied pegs in the corresponding configuration.
          This leads to a natural decomposition of the vertex set into
          \textit{occupancy classes}.
        \item We provide a detailed combinatorial decomposition of the edge set of \(H_p^n\) according to occupancy transitions between pegs, leading to an exact counting of all degree pairs.
    \item We express all coefficients of \(M(H_p^n; x, y)\) in terms of Stirling numbers of the second kind and falling factorials, ensuring computational tractability for arbitrary \(p\) and \(n\).

    \item Combining these ingredients, we  derive a
          fully explicit closed formula for $M(H_p^n;x,y)$.
    \item We verify the correctness of our formulas through explicit enumeration for small instances and demonstrate how classical degree-based indices (Zagreb, Randić, etc.) can be directly obtained via differential operators.
\end{itemize}

The remainder of the paper is organized as follows. 
Section~\ref{sec:preliminaries} recalls basic definitions and the operator 
formulae linking the M-polynomial to the most common degree-based indices.
Section~\ref{sec:hanoi} introduces generalized Hanoi graphs and the 
occupancy-based viewpoint.
Section~\ref{sec:analysis} develops the degree-based combinatorial
framework, culminating in a full decomposition of the edge set.
Section~\ref{sec:explicit} presents the explicit formula for 
$M(H_p^n;x,y)$.
Section~\ref{sec:numerical} provides numerical validation and illustrative 
examples.  
Finally, Section~\ref{sec:conclusion} discusses perspectives and possible 
extensions, including restricted Hanoi graphs and other self-similar families.

\section{Preliminaries}\label{sec:preliminaries}

Throughout this paper, all graphs are assumed to be finite, simple, and connected.  
For a graph \(G = (V,E)\), we denote by \(V(G)\) its vertex set and by \(E(G)\) its edge set.  
An edge \(e = \{u,v\} \in E(G)\) connects vertices \(u,v \in V(G)\), and the \emph{degree} of a vertex \(u\) is denoted by \(\deg(u)\), i.e.,
\[
\deg(u) = |\{ v \in V(G) : \{u,v\} \in E(G)\}|.
\]
The minimum and maximum degrees of \(G\) are denoted respectively by \(\delta(G)\) and \(\Delta(G)\).

A \emph{molecular graph} is a graph-theoretical abstraction of a chemical compound, in which vertices correspond to non-hydrogen atoms and edges represent covalent bonds \cite{gutman2012mathematical}.  
Topological indices derived from such graphs provide quantitative descriptors of molecular structure and are used in quantitative structure–property (QSPR) and structure–activity (QSAR) studies \cite{Dearden2017}.

A \emph{graph invariant} is a function \(I\) that assigns a unique value \(I(G)\) to each graph \(G\) such that isomorphic graphs have equal values.  
When this invariant depends only on the vertex degrees, it is called a \emph{degree-based topological index}.

\begin{definition}[M-polynomial {\cite{EmericMploy}}]
The \emph{M-polynomial} of a graph \(G\) is defined as
\begin{equation}
M(G; x, y)
 = \sum_{i \le j} m_{i,j}(G) x^i y^j,
\end{equation}
where \(m_{i,j}(G)\) denotes the number of edges \(\{u,v\} \in E(G)\) such that
\(\{\deg(u), \deg(v)\} = \{i,j\}\); that is,
\[
m_{i,j}(G)
= \big| \{\, \{u,v\} \in E(G) \mid \{\deg(u),\deg(v)\} = \{i,j\} \,\} \big|.
\]
\end{definition}
The strength of the M-polynomial lies in its universality: once \(M(G; x,y)\) is known,  
a broad class of classical indices can be obtained by applying differential or integral operators to it.  
Table~\ref{tab:indices} summarizes the most common degree-based indices and their operator expressions.

For completeness, we recall several widely used degree-based topological indices  that can be derived from the M-polynomial:

\begin{definition}
The \emph{first Zagreb}, \emph{second Zagreb}, and \emph{second modified Zagreb} indices of $G$ are defined respectively as:
\begin{align}
    M_1(G) &= \sum_{\{u,v\} \in E(G)} \left(\deg(u) + \deg(v)\right), \\
    M_2(G) &= \sum_{\{u,v\} \in E(G)} \deg(u) \cdot \deg(v), \\
    MM_2(G) &= \sum_{\{u,v\} \in E(G)} \frac{1}{\deg(u) \cdot \deg(v)}.
\end{align}
\end{definition}

\begin{definition}
The \emph{generalized Randić} and \emph{reciprocal generalized Randić} indices of $G$ are defined, respectively, as:
\begin{align}
    R_\alpha(G) &= \sum_{\{u,v\} \in E(G)} \left(\deg(u) \cdot \deg(v)\right)^{\alpha}, \\
    RR_\alpha(G) &= \sum_{\{u,v\} \in E(G)} \frac{1}{\left(\deg(u) \cdot \deg(v)\right)^{\alpha}}.
\end{align}
\end{definition}

\begin{definition}
The \emph{symmetric division index} of $G$ is defined by:
\begin{equation}
    SSD(G) = \sum_{\{u,v\} \in E(G)} \left( \frac{\min\{\deg(u), \deg(v)\}}{\max\{\deg(u), \deg(v)\}} + \frac{\max\{\deg(u), \deg(v)\}}{\min\{\deg(u), \deg(v)\}} \right).
\end{equation}
\end{definition}

\begin{definition}
The \emph{harmonic}, \emph{inverse sum indeg}, \emph{augmented Zagreb}, and \emph{Forgotten} indices of $G$ are defined respectively as:
\begin{align}
    H(G) &= \sum_{\{u,v\} \in E(G)} \frac{2}{\deg(u) + \deg(v)}, \\
    ISI(G) &= \sum_{\{u,v\} \in E(G)} \frac{\deg(u) \cdot \deg(v)}{\deg(u) + \deg(v)}, \\
    A(G) &= \sum_{\{u,v\} \in E(G)} \left( \frac{\deg(u) \cdot \deg(v)}{\deg(u) + \deg(v) - 2} \right)^3, \\
    F(G) &= \sum_{\{u,v\} \in E(G)} \left( \deg(u)^2 + \deg(v)^2 \right).
\end{align}
\end{definition}

Let \(D_x, D_y, S_x, S_y, J, Q_\alpha\) denote the operators defined by

\begin{align*}
    D_x(g(x,y)) &= x \frac{\partial g(x,y)}{\partial x}, & D_y(g(x,y)) &= y \frac{\partial g(x,y)}{\partial y}, \\
    S_x(g(x,y)) &= \int_{0}^{x} \frac{g(t,y)}{t} \, dt, & S_y(g(x,y)) &= \int_{0}^{y} \frac{g(x,t)}{t} \, dt, \\
    J(g(x,y)) &= g(x,x), & Q_\alpha(g(x,y)) &= x^{\alpha} g(x,y).
\end{align*}
Then the indices above satisfy the operator identities in Table~\ref{tab:indices}.

\begin{table}[htbp]
\centering

\begin{tabular}{ll}
\toprule
Index & Relation with \(M(G;x,y)\) \\
\midrule
First Zagreb & \(M_1(G) = (D_x + D_y) M(G;x,y)\big|_{x=y=1}\) \\
Second Zagreb & \(M_2(G) = (D_x D_y) M(G;x,y)\big|_{x=y=1}\) \\
Modified second Zagreb & \(MM_2(G) = (S_x S_y) M(G;x,y)\big|_{x=y=1}\) \\
Generalized Randić & \(R_{\alpha}(G) = (D_x^\alpha D_y^\alpha) M(G;x,y)\big|_{x=y=1}\) \\
Reciprocal generalized Randić & \(RR_{\alpha}(G) = (S_x^\alpha S_y^\alpha) M(G;x,y)\big|_{x=y=1}\) \\
Symmetric division & \(SSD(G) = (D_x S_y + D_y S_x) M(G;x,y)\big|_{x=y=1}\) \\
Harmonic & \(H(G) = 2 S_x J(M(G;x,y))\big|_{x=1}\) \\
Inverse sum indeg & \(ISI(G) = S_x J D_x D_y M(G;x,y)\big|_{x=1}\) \\
Augmented Zagreb & \(A(G) = S_x^3 Q_{-2} J D_x^3 D_y^3 M(G;x,y)\big|_{x=1}\) \\
Forgotten & \(F(G) = (D_x^2 + D_y^2) M(G;x,y)\big|_{x=y=1}\) \\
\bottomrule
\end{tabular}
\caption{Relations between degree-based indices and the M-polynomial.}
\label{tab:indices}
\end{table}

\section{Generalized Hanoi graphs}\label{sec:hanoi}

The Tower of Hanoi (TH) problem, invented by the French number theorist Édouard Lucas in 1883, has long served as a source of challenge and inspiration in mathematics, computer science. For more details on the mathematics related to the Tower of Hanoi, we refer the reader to the comprehensive monograph by Hinz \emph{et al.}~\cite{hinz2018tower}.

In its general form, the Tower of Hanoi problem is defined for $p \in \mathbb{N}$, with $p \geq 3$ pegs and $n \in \mathbb{N}_0$ discs, each of a distinct size. A  \textit{legal move} is a transfer of the topmost disc\footnote{A topmost disc is a top disc of a peg.}  from one peg to another, without ever placing a larger disc on top of a smaller one. Initially, all discs are stacked on a single peg in decreasing order of size from bottom to top—this is called a \textit{perfect state}. The objective is to move all discs from one perfect state to another using the minimum number of legal moves.

A state  (i.e., a distribution of discs across the pegs) is called   \textit{regular}  if the discs on each peg are arranged in decreasing size from bottom to top. The Hanoi graphs $H_p^n$ provide a natural mathematical model for the TH problem. In these graphs, each vertex represents a regular state, and two vertices are connected by an edge if one can be obtained from the other through a single legal move.

Numerous structural and combinatorial properties of Hanoi graphs have been explored in~\cite{hinz2018tower} and the references therein.

\begin{definition}[\cite{hinz2018tower}]\label{def:hanoi_graph_defenition}
The Hanoi graphs $H_p^n$ for base $p \in \mathbb{N}_3$ and exponent $n \in \mathbb{N}_0$ are defined as follows.
    \begin{align*}
        V(H_p^n)&=\Set{s_{n}\cdots s_{1}\mid s_{d}\in[p]_{0},\;  d\in [n]}\\
        E(H_p^n)&=\Set{\{ \underline{s}i\overline{s},\underline{s}j\overline{s}\}\mid i,j\in[p]_{0},i\neq j, \underline{s}\in [p]_{0}^{n-d},\overline{s}\in ([p]_{0}\setminus\{i,j\})^{d-1}, d\in [n]  }
    \end{align*}
\end{definition}

\input{small_graphs}

\begin{definition}[Occupancy and structural parameters]
    For a state \(s \in V(H_p^n)\), define the occupancy number $o(s) : V(H_p^n) \rightarrow [p]_0$ by
\[
o(s) = |\{ i \in [p] \;|\; \exists\, d \in [n] : s_d = i \}|,
\]
i.e., the number of pegs that  contain at least one disc in state $s$.  
\end{definition}

Clearly, if $n = 0$, then $V(H_p^n) = \emptyset$, and the function $o(s)$ is undefined (or vacuously equal to zero). Otherwise, 

\[
1 \le o(s) \le r = \min\{n,p\}.
\]

\begin{table}[H]
\centering

\renewcommand{\arraystretch}{1.25}
 \resizebox{\columnwidth}{!}{
\begin{tabular}{ll}
\toprule
\textbf{Symbol} & \textbf{Meaning} \\
\midrule
$n$ 
& Number of discs \\
$p$ 
& Number of pegs \\
$r = \min\{n, p\}$ & Maximum possible number of occupied pegs\\

\midrule
$H_p^n$ 
& Generalized Hanoi graph with $p$ pegs and $n$ discs \\

$V(H_p^n)$, $E(H_p^n)$ 
& Vertex set and edge set of $H_p^n$ \\

$o(s)$ 
& Occupancy of state $s$ (number of pegs containing at least one disc) \\

$\Omega_p^n = \{1,\dots,r\}$ 
& Set of possible occupancy values, where $r = \min\{n,p\}$ \\

$D_\mu = \{s \in V(H_p^n) : o(s)=\mu\}$ 
& Class of configurations with exactly $\mu$ occupied pegs \\

$f_{p,n}(\mu)$ 
& Degree of any vertex with occupancy $\mu$ \\

$\lambda_\mu = f_{p,n}(\mu)$ 
& Degree value associated with occupancy $\mu$ \\

$\delta(H_p^n)$, $\Delta(H_p^n)$ 
& Minimum and maximum vertex degrees of $H_p^n$ \\

\midrule

$O_p^n(\mu)$ 
& Number of configurations with $\mu$ occupied pegs \\

$O_p^n(\nu \mid \mu)$ 
& Number of configurations with $\mu$ occupied pegs, $\nu$ of which are singletons \\

$p^{\underline{\mu}} = p(p-1)\cdots(p-\mu+1)$ 
& Falling factorial of order $\mu$ \\

$\begin{Bmatrix} n \\ k \end{Bmatrix}$ 
& Stirling number of the second kind \\

$\begin{Bmatrix} n \\ k \end{Bmatrix}_{\ge 2}$ 
& 2-associated Stirling number (each block has size $\ge 2$) \\

\midrule

$A_1(H_p^n)$ 
& Edges corresponding to moves OP $\to$ OP \\

$A_2(H_p^n)$ 
& Edges corresponding to moves OP $\to$ EP \\

$E_1(H_p^n)$ 
& Edges from singleton peg $\to$ empty peg \\

$E_2(H_p^n)$ 
& Edges from multiton peg $\to$ occupied peg \\

$E_3(H_p^n)$ 
& Edges from multiton peg $\to$ empty peg (increasing occupancy) \\

\midrule

$M(H_p^n; x, y)$ 
& M-polynomial of $H_p^n$ \\

$m_{\lambda,\lambda'}(H_p^n)$ 
& Number of edges joining degree classes $\lambda$ and $\lambda'$ \\

$\Gamma_p^n = \{\lambda_\mu : \mu \in \Omega_p^n\}$ 
& Set of all realized vertex degrees in $H_p^n$ \\

\bottomrule
\end{tabular}
}
\caption{Summary of notation used throughout the paper.}
\label{tab:notation}
\end{table}

 We denote by \(\Omega_p^n = [r]\) the set of possible occupancy values.  Pegs may therefore be classified as follows:
\begin{itemize}
    \item \textit{Empty peg } (EP): contains no disc.
    \item \textit{Occupied Peg} (OP): contains at least one disc.
    \item \textit{Singleton peg } (SP): is an OP that contains exactly one disc.
    \item \textit{Multiton peg } (MP): is an OP that contains two or more discs.
   
\end{itemize}

Table~\ref{tab:notation} summarizes the notation used throughout the paper for ease of reference.

\section{Degree-based combinatorial analysis}\label{sec:analysis}

In this section, we develop a detailed combinatorial characterization of vertex degrees and edge structures in the generalized Hanoi graphs \(H_p^n\).  
We begin by analyzing how vertex degrees depend on the occupancy number of the configuration, then enumerate all vertices according to their occupancy composition, and finally derive explicit formulas for the number of edges of different structural types.  
This analysis constitutes the combinatorial foundation of the M-polynomial derivation that follows.

\subsection{Dependence of vertex degree on occupancy}

\begin{lemma}[Degree as a function of occupancy \cite{HINZ20121521,Arett01062010}]
\label{lem:deg_occ}
For a vertex \(s\) with occupancy \(o(s)=\mu\), its degree in \(H_p^n\) is given by
\begin{equation}
\label{eq:f_pn}
\deg(s) = f_{p,n}(\mu)
  = \sum_{\ell=1}^{\mu}(p-\ell)
  = \binom{p}{2}-\binom{p-\mu}{2}.
\end{equation}
\end{lemma}

\begin{proof}
Consider a configuration \(s\) in which \(\mu\) pegs are occupied.  
Each of the \(\mu\) pegs contains a stack of discs, and the smallest (topmost) disc on each peg can be moved to another peg, provided that the destination peg either is empty or has a larger top disc.

Let us enumerate all possible target pegs for each of these \(\mu\) possible moves.  
For the peg containing the smallest top disc overall, all the remaining \(p-1\) pegs are available destinations.  
For the peg containing the next smallest top disc, one of these destinations is forbidden (the peg holding the globally smallest top disc), leaving \(p-2\) options.  
Proceeding in this manner, the total number of legal moves (hence, the degree of \(s\)) equals
\[
(p-1)+(p-2)+\cdots+(p-\mu)
 = \sum_{\ell=1}^{\mu}(p-\ell)
 = \binom{p}{2}-\binom{p-\mu}{2}.
\]
This argument relies on the symmetry and legality constraints of the Tower of Hanoi rules and holds for any valid configuration with occupancy \(\mu\).\qedhere
\end{proof}

\begin{corollary}[Degree bounds]
\[
\delta(H_p^n)=f_{p,n}(1)=p-1, \qquad
\Delta(H_p^n)=f_{p,n}(r)=\binom{p}{2}-\binom{p-r}{2}.
\]
\end{corollary}

\begin{proof}
The smallest degree arises when only one peg is occupied (\(\mu=1\)), since only that single peg provides possible moves to the other \(p-1\) pegs.  
Conversely, the largest degree occurs when all \(r=\min\{n,p\}\) pegs are occupied, producing the maximum number of active move pairs.\qedhere
\end{proof}

\begin{corollary}[Characterization of the degree set of $H_p^n$]\label{lemma:degree-set}
For all $p \in \mathbb{N}_3$ and $n \in \mathbb{N}_0$, the set of all possible vertex degrees in the graph $H_p^n$ is given by
\begin{equation}
    \Gamma_p^n 
    = \left\{ \sum_{\ell = 1}^{\mu} (p - \ell) \;\middle|\; \mu \in \Omega_p^n \right\}
    = \left\{\, p - 1,\; 2p - 3,\; \ldots,\; \sum_{\ell = 1}^{\min\{n, p\}} (p - \ell) \,\right\}.
\end{equation}
\end{corollary}

\begin{proof}
The result follows directly from the proof of Lemma~\ref{lem:deg_occ}.\qedhere
\end{proof}

\begin{proposition}[Non-injectivity and surjectivity of the mapping $f_{p,n}$]\label{prop:deg-noninj}
For all $p \in \mathbb{N}_3$ and $n \in \mathbb{N}_0$, the function 
\[
f_{p,n} : \Omega_p^n \longrightarrow \Gamma_p^n
\]
is surjective but not injective. Consequently, $f_{p,n}$ is not bijective.
\end{proposition}

\begin{proof}
From Lemma~\ref{lemma:Deltarepeated}, we know that
\[
f_{p,n}(p) = f_{p,n}(p - 1) = \Delta(H_p^n),
\]
while clearly $p \neq p - 1$. Hence, $f_{p,n}$ maps two distinct elements of $\Omega_p^n$ to the same value, and is therefore not injective.

To show surjectivity, let $\gamma \in \Gamma_p^n$ and consider the equation $f_{p,n}(\mu) = \gamma$.  

\begin{itemize}
    \item If $\gamma \in \Gamma_p^n \setminus \{\Delta(H_p^n)\}$, then
    \[
    \mu = p - \frac{1 + \sqrt{\,1 + 4\bigl(p(p - 1) - 2\gamma\bigr)\,}}{2}
    \]
    yields a unique $\mu \in \Omega_p^n$ such that $f_{p,n}(\mu) = \gamma$.
    
    \item If $\gamma = \Delta(H_p^n)$, there exist exactly two values $\mu_1, \mu_2 \in \Omega_p^n$ satisfying
    \[
    \mu_1 = p - \frac{1 + \sqrt{\,1 + 4\bigl(p(p - 1) - 2\gamma\bigr)\,}}{2}, 
    \qquad
    \mu_2 = p - \frac{1 - \sqrt{\,1 + 4\bigl(p(p - 1) - 2\gamma\bigr)\,}}{2},
    \]
    and both verify $f_{p,n}(\mu_1) = f_{p,n}(\mu_2) = \gamma$.
\end{itemize}

Thus, every $\gamma \in \Gamma_p^n$ admits at least one preimage in $\Omega_p^n$, which proves that $f_{p,n}$ is surjective.\qedhere
\end{proof}

\begin{remark}
The degree of a vertex depends exclusively on its occupancy \(\mu=o(s)\); the exact arrangement of discs among the \(\mu\) pegs is irrelevant.  
Therefore, all vertices can be partitioned into \emph{degree classes}
\[
D_\mu=\{\,s\in V(H_p^n)\mid o(s)=\mu\,\},
\]
and each vertex in \(D_\mu\) has degree \(f_{p,n}(\mu)\).
\end{remark}

\subsection{Adjacency relations between degree classes}

We next determine which degree classes are connected by edges in \(H_p^n\).

\begin{proposition}[Adjacency of occupancy classes]
\label{prop:adj_occ}
If \(u,v\in V(H_p^n)\) are adjacent, then
\[
|o(u)-o(v)|\in\{0,1\}.
\]
\end{proposition}

\begin{proof}
A legal move transfers one disc from a source peg to a target peg.

\emph{Case 1.}  
If the target peg already contains one or more discs, the total number of occupied pegs does not increase.  
If the source peg still contains other discs after the move, the total occupancy remains unchanged; if it becomes empty, the number of occupied pegs decreases by one.  

\emph{Case 2.}  
If the target peg was empty, then after the move it becomes occupied, possibly increasing the number of occupied pegs by one (unless the source peg was a singleton, in which case the occupancy remains constant).

Hence the occupancy can only remain the same or differ by one between adjacent vertices.\qedhere
\end{proof}

\begin{corollary}[Difference of degrees between adjacent vertices]\label{cor:deg-step}
If \(u,v\in V(H_p^n)\) are adjacent and \(o(u)\le o(v)\), then

\begin{equation}
    \deg(u)=
\begin{cases}
\deg(v), & \text{if } o(u)=o(v),\\[4pt]
\deg(v)-(p-o(v)), & \text{if } o(v)=o(u)+1.
\end{cases}
\end{equation}

\end{corollary}

\begin{proof}
When \(o(u)=o(v)\), both vertices belong to the same degree class and therefore share the same degree \(f_{p,n}(o(u))\).  
If \(o(v)=o(u)+1\), then from \eqref{eq:f_pn},
\[
\deg(v)-\deg(u)=f_{p,n}(o(u)+1)-f_{p,n}(o(u))
=(p-o(v)),
\]
yielding the claimed difference.\qedhere
\end{proof}

\begin{corollary}
If \(u, v \in V(H_p^n)\) are adjacent and either \(o(u) = p\) or \(o(v) = p\), then \(\deg(u) = \deg(v)\).
\end{corollary}

\begin{proof}
The result follows directly from Corollary~\ref{cor:deg-step}.
\end{proof}

\begin{lemma}[Characterization of vertices of maximal degree]\label{lemma:Deltarepeated}
For all $p \in \mathbb{N}_3$, $n \in \mathbb{N}_0$, and $s \in V(H_p^n)$, we have

\begin{equation}
    \deg(s) = \Delta(H_p^n)
\quad \Longleftrightarrow \quad
o(s) \in \{p,\, p-1\}.
\end{equation}

\end{lemma}

\begin{lemma}[Monotonicity between occupancy and degree]\label{lemma:O_D}
For all $p \in \mathbb{N}_3$, $n \in \mathbb{N}_0$, and $\{u, v\} \in E(H_p^n)$, the following implication holds:
\begin{equation}
    o(u) \leq o(v) \ \Longrightarrow\ \deg(u) \leq \deg(v).
\end{equation}
\end{lemma}

\begin{proof}
If $o(u) = o(v)$, then by Lemma~\ref{lem:deg_occ}, we have $\deg(u) = \deg(v)$, since the degree of a state depends solely on the number of occupied pegs.

Now assume, without loss of generality, that $o(u) = o(v) - 1$ and $o(u) < p$.  
Let the $o(v)$ topmost discs in state $v$ be indexed in increasing order of size as $d_1, d_2, \dots, d_{o(v)}$.  
Consider removing the largest disc $d_{o(v)}$; the remaining $o(v) - 1 = o(u)$ discs then generate exactly $\deg(u)$ legal moves.  
The additional disc $d_{o(v)}$ contributes $p - o(v)$ further legal moves, since it can move to all pegs not already occupied by one of the smaller discs.

Hence,
\[
\deg(v) = \deg(u) + (p - o(v)) \geq \deg(u),
\]
which proves the claim.
\end{proof}

\begin{proposition}[Difference of degrees between adjacent vertices]\label{prop:deg-diff}
Let $p \in \mathbb{N}_3$, $n \in \mathbb{N}_0$, and $\{u, v\} \in E(H_p^n)$. Then
\begin{equation}
    |\deg(u) - \deg(v)| \in \left\{\, 0,\; p - \max\{o(u),\, o(v)\} \,\right\}.
\end{equation}
\end{proposition}

\begin{proof}
If \(o(u)=o(v)\), the difference is zero.  
If \(o(v)=o(u)+1\), then from the proof of Lemma~\ref{lemma:O_D},
\(\deg(v)-\deg(u)=p-o(v)\), which equals \(p-\max\{o(u),o(v)\}\).  
The same holds symmetrically when \(o(u)=o(v)+1\).
\end{proof}

\subsection{Enumeration of configurations by occupancy}

Let \(O_p^n(\mu)\) denote the number of vertices in \(H_p^n\) having exactly \(\mu\) occupied pegs.  
We derive a closed formula for this quantity.

\begin{theorem}[Number of configurations by occupancy \cite{klavzarcombinatorics,Arett01062010}]
\label{th:Opmn}
For all \(p\ge3\) and \(n\ge0\),
\begin{equation}
\label{eq:Opmn}
O_p^n(\mu)=\begin{Bmatrix} n \\ \mu \end{Bmatrix}\,p^{\underline{\mu}},
\qquad 1\le\mu\le r,
\end{equation}
where \(\begin{Bmatrix} n \\ \mu \end{Bmatrix}\) is the Stirling number of the second kind and \(p^{\underline{\mu}}=p(p-1)\cdots(p-\mu+1)\) the falling factorial.
\end{theorem}

\begin{proof}
We must count all distinct distributions of \(n\) labeled discs among \(p\) labeled pegs such that exactly \(\mu\) pegs are nonempty.

First, we partition the \(n\) discs into \(\mu\) nonempty groups, each corresponding to a stack on one peg.  
The number of such partitions is \(\begin{Bmatrix} n \\ \mu \end{Bmatrix}\).  
Then, we assign each group to a distinct peg.  
Since the pegs are labeled, this can be done in \(p^{\underline{\mu}}=p!/(p-\mu)!\) ways.

Each assignment yields a distinct configuration because the order of discs within a stack is determined by their sizes.  
Multiplying both factors gives the result \eqref{eq:Opmn}.
\end{proof}

\begin{lemma}[\cite{klavzarcombinatorics}]\label{lem:sum-O=pn}
For all \(p\ge3\) and \(n\ge0\), let $r = \min\{n,p\}$, then  

\begin{equation}
    \sum_{\mu=1}^{r} O_p^n(\mu) = p^n.
\end{equation}
\end{lemma}

\begin{proof}
This follows from the classical identity
\(\sum_{\mu=0}^n \begin{Bmatrix} n \\ \mu \end{Bmatrix}x^{\underline{\mu}} = x^n\)
with \(x=p\).
\end{proof}

We refine the counting by distinguishing the number of singleton pegs within a configuration. Let \(O_p^n(\nu\mid\mu)\) denote the number of configurations with \(\mu\) occupied pegs, among which \(\nu\) pegs are singletons (contain exactly one disc).

\begin{theorem}[Configurations with \(\mu\) occupied and \(\nu\) singletons]
\label{th:refined}
For all \(1\le\nu\le\mu\le r\),
\begin{equation}
\label{eq:refined}
O_p^n(\nu\mid\mu)
=\binom{n}{\nu}\,
  \begin{Bmatrix} n-\nu\\ \mu-\nu \end{Bmatrix}_{\geq 2}\,
  p^{\underline{\mu}},
\end{equation}
where \(\begin{Bmatrix} n-\nu\\ \mu-\nu \end{Bmatrix}_{\geq 2}\) denotes the 2-associated Stirling number of the second kind.
\end{theorem}

\begin{proof}
We first choose \(\nu\) discs to serve as singletons, in \(\binom{n}{\nu}\) ways.  
The remaining \(n-\nu\) discs must be arranged on the remaining \(\mu-\nu\) pegs, each receiving at least two discs.  
The number of such partitions is \(\begin{Bmatrix} n-\nu\\ \mu-\nu \end{Bmatrix}_{\geq 2}\).  
Finally, we assign these \(\mu\) stacks to distinct pegs in \(p^{\underline{\mu}}\) ways.  
The product of these three factors yields \eqref{eq:refined}.
\end{proof}

\begin{lemma}
\label{cor:split-Opnu}
Summing over all admissible values of \(\nu\),

\begin{equation}
O_p^n(\mu)=\sum_{\nu=0}^{\mu}O_p^n(\nu\mid\mu)
=\sum_{\nu=0}^{\mu}\binom{n}{\nu}
   \begin{Bmatrix} n-\nu\\ \mu-\nu \end{Bmatrix}_{\geq 2}\,p^{\underline{\mu}}.
\end{equation}

with the convention that
\[
\begin{Bmatrix} n - \mu \\ 0 \end{Bmatrix}_{\geq 2} =
\begin{cases}
1, & \text{if } n = \mu; \\
0, & \text{otherwise}.
\end{cases}
\]
\end{lemma}

\subsection{Edge decomposition and degree-based classification}

In this section, we refine the combinatorial structure of the edge set of the generalized Hanoi graph $H_p^n$ by decomposing it into several disjoint subsets according to the nature of the moves between states.  
We begin by recalling that each edge $\{u,v\} \in E(H_p^n)$ corresponds to a legal move of a disc between two pegs under the standard Tower of Hanoi rules.  
Depending on whether this move increases, decreases, or preserves the number of occupied pegs, the edges can be classified into different categories.

\noindent
The following Proposition establishes a first-level partition of the edge set into two fundamental blocks.

\begin{proposition}\label{prop:two-block}
For all $p \in \mathbb{N}_3$ and $n \in \mathbb{N}_0$, the edge set of the Hanoi graph $H_p^n$ can be expressed as
\begin{equation}
    E(H_p^n) = A_1(H_p^n) \cup A_2(H_p^n),
\end{equation}
where $A_1(H_p^n)\cap  A_2(H_p^n)=\emptyset$, and
\begin{align}
     A_1(H_p^n)&=\left\{\{u,v\}\in E(H_p^n) \mid \text{\{u,v\} corresponds to a move from an OP to another OP} \right\},\\
     A_2(H_p^n)&=\left\{\{u,v\}\in E(H_p^n) \mid \text{\{u,v\} corresponds to a move from an OP to an EP} \right\}.
\end{align}
\end{proposition}

The next result provides closed expressions for the cardinalities of these two classes.

\begin{proposition}\label{prop:AA}
For all $p \in \mathbb{N}_3$ and $n \in \mathbb{N}_0$,  let $r = \min\{n, p\}$. Then we have
\begin{align}
     |A_1(H_p^n)|&=\dfrac{1}{2}\sum_{\mu=1}^{r}\binom{\mu}{2}O_p^n(\mu),\\
     |A_2(H_p^n)|&=\dfrac{1}{2}\sum_{\mu=1}^{r}\mu(p-\mu)O_p^n(\mu).
\end{align}
\end{proposition}

\noindent
From these two disjoint families, one can recover the total number of edges of $H_p^n$.  
This value has been previously obtained by Klavžar and Milutinović~\cite{klavzarcombinatorics} and Arett~\cite{Arett01062010} through a different reasoning based on the degree distribution of the vertices.  
Their approach consisted in summing, for each possible number of occupied pegs $\mu$, the product of the number of corresponding states $O_p^n(\mu)$ and the degree $f_{p,n}(\mu)$ of each state, giving:
\begin{equation}\label{eq:Ehpn-sum-mu}
    E(H_p^n)=\dfrac{1}{2} \sum_{\mu=1}^{r}O_{p}^{n}(\mu)f_{p,n}(\mu)
    =\dfrac{1}{2}\sum_{\mu=1}^{r}\begin{Bmatrix} n \\ \mu\end{Bmatrix}p^{\underline{\mu}}\left[\binom{p}{2}-\binom{p-\mu}{2}\right].
\end{equation}
A similar method was later extended to restricted Hanoi graphs in~\cite{mehiri2024onrestricted}, a result that can also be applied to the classical family $H_p^n$.

\noindent
However, by exploiting the disjoint decomposition introduced in Proposition~\ref{prop:two-block}, we can recover the same result more structurally, as shown below.

\begin{theorem}[\cite{klavzarcombinatorics,Arett01062010}]\label{theo:edges-sum-forms}  
For all $p \in \mathbb{N}_3$ and $n \in \mathbb{N}_0$, let $r = \min\{n, p\}$. Then we have
\begin{equation}\label{eq:Ehpn-alt}
        E(H_p^n)=\dfrac{1}{4}\sum_{\mu=1}^{r}\mu\left(2p-\mu-1 \right) O_p^n(\mu).
\end{equation}
\end{theorem}
\begin{proof}
From Propositions~\ref{prop:two-block} and~\ref{prop:AA}, we have 
\begin{align*}
E(H_p^n)&=|A_1(H_p^n)|+|A_2(H_p^n)|\\
&=\dfrac{1}{2}\sum_{\mu=1}^{r}\binom{\mu}{2}O_p^n(\mu)+\dfrac{1}{2}\sum_{\mu=1}^{r}\mu(p-\mu)O_p^n(\mu)\\
&=\dfrac{1}{2}\sum_{\mu=1}^{r}\left(\dfrac{\mu(\mu-1)}{2}+\mu(p-\mu) \right) O_p^n(\mu) \\
&=\dfrac{1}{4}\sum_{\mu=1}^{r}\mu\left(2p-\mu-1 \right) O_p^n(\mu).\qedhere
\end{align*}
\end{proof}

We now refine the above dichotomy by examining more closely the nature of the move that generates each edge.  
Depending on whether the source and target configurations differ by one in the number of occupied pegs, or whether the moved disc comes from a singleton or multiton peg, we distinguish three mutually exclusive situations.
\begin{theorem}\label{thm:edge_partition}
Let $p \in \mathbb{N}_3$ and $n \in \mathbb{N}_0$.  
For every edge $\{u, v\} \in E(H_p^n)$ with $o(u) \leq o(v)$. Then the edge $\{u, v\}$ corresponds to exactly one of the following three disjoint cases:
\begin{itemize}
    \item[$(\texttt{c}_{1})$] A disc is moved from a singleton peg to an empty peg (occurs only if $o(u) = o(v)$);
    \item[$(\texttt{c}_{2})$] A disc is moved from a multiton peg to an occupied peg (also only if $o(u) = o(v)$);
    \item[$(\texttt{c}_{3})$] A disc is moved from a multiton peg to an empty peg (occurs only if $o(u) = o(v) - 1$).
\end{itemize}

Hence, the edge set of $H_p^n$ admits the following disjoint decomposition:
\begin{equation}
    E(H_p^n) = E_1(H_p^n) \cup E_2(H_p^n) \cup E_3(H_p^n),
\end{equation}
where
\begin{align}
    E_i(H_p^n) &= \{\{u, v\} \in E(H_p^n) \mid \text{$\{u,v\}$ satisfies \texttt{(c\textsubscript{$i$})}}\}, \quad i = 1, 2, 3,
\end{align}
and these subsets are pairwise disjoint, i.e., $E_i(H_p^n) \cap E_j(H_p^n) = \emptyset$ for all $i \neq j$.
\end{theorem}

\begin{proof}
This result follows directly from the case analysis presented in the proof of Proposition~\ref{prop:adj_occ}.
\end{proof}

Note that the reverse situations of \texttt{(c\textsubscript{2})} and \texttt{(c\textsubscript{3})} cannot occur under the assumption \(o(u) \leq o(v)\), as they would imply a decrease in the number of occupied pegs. Consequently, only case \texttt{(c\textsubscript{1})} is symmetric. This observation will be taken into account in ulterior enumerations.

The following theorem enumerates the number of edges belonging to each of these three categories.

\begin{theorem}
For all $p \in \mathbb{N}_3$ and $n \in \mathbb{N}_0$,  let $r = \min\{n, p\}$. Then we have
\begin{align}
    |E_1(H_p^n)| &= \dfrac{1}{2} \sum_{\mu = 1}^{r} (p - \mu) \sum_{\nu = 0}^{\mu} \nu \, O_p^n(\nu \mid \mu),\label{eq:E1} \\
    |E_2(H_p^n)| &= \dfrac{1}{2} \binom{p}{2} p^n - \dfrac{1}{4} \sum_{\mu = 1}^{r} (p - \mu) \sum_{\nu = 0}^{\mu} \left(p + 3\mu - 2\nu - 1 \right) O_p^n(\nu \mid \mu),\label{eq:E2} \\
    |E_3(H_p^n)| &= \sum_{\mu = 1}^{r} (p - \mu) \sum_{\nu = 0}^{\mu} (\mu - \nu) \, O_p^n(\nu \mid \mu).\label{eq:E3}
\end{align}
\end{theorem}

\begin{proof}
The term $|E_1(H_p^n)|$ counts edges $\{u, v\} \in E(H_p^n)$ such that $o(u) = o(v) \leq p - 1$, where the move is from a singleton peg to an empty peg.  
The number of such states $u$ is given by $O_p^n(\nu \mid \mu)$, where $\mu \in [p - 1]$ and $\nu \in [\mu]$.  
For each such state $u$ (or equivalently $v$), there are $\nu$ singleton pegs from which a disc can be chosen, and $p - \mu$ empty pegs on which it can be placed.  
Thus, the number of such moves is $\nu (p - \mu)$.  
Since edges are unordered, we divide by $2$ to remove symmetry, which yields Equation~\eqref{eq:E1}.

The term $|E_3(H_p^n)|$ counts edges $\{u, v\} \in E(H_p^n)$ where $o(v) = o(u) + 1 \leq p$, and the move is from a multiton peg to an empty peg in state $u$.  
Such a move is only possible when $\mu \leq p - 1$ and $\nu < \mu$, since having all occupied pegs as singletons ($\nu = \mu$) would exclude the presence of multitons.  
For each such state $u$, there are $\mu - \nu$ multitons from which a disc can be chosen, and $p - \mu$ empty pegs available to receive it.  
Moreover, no symmetry elimination is required since the directionality is fixed by $o(u) \neq o(v)$.  
This gives Equation~\eqref{eq:E3}.

The term $|E_2(H_p^n)|$ accounts for edges with $o(u) = o(v)$ where the move is from a multiton peg to another occupied peg.  
A direct enumeration of these configurations is infeasible due to dependencies on disc sizes at both the source and target pegs.  
Instead, we apply the inclusion–exclusion principle.  
From Theorem~\ref{thm:edge_partition}, we have
\[
|E_2(H_p^n)| = |E(H_p^n)| - |E_1(H_p^n)| - |E_3(H_p^n)|.
\]
Substituting the expression of $|E(H_p^n)|$ from Theorem~\ref{theo:edges-sum-forms}, together with those of $|E_1(H_p^n)|$ and $|E_3(H_p^n)|$, yields
\begin{align*}
    |E_2(H_p^n)| &= \frac{1}{2} \sum_{\mu = 1}^{r} \sum_{\nu = 0}^{\mu} 
    O_p^n(\nu \mid \mu) \!\left( \binom{p}{2} - \binom{p - \mu}{2} \right) \\
    &\quad - \frac{1}{2} \sum_{\mu = 1}^{r} \sum_{\nu = 0}^{\mu} 
    \nu (p - \mu)\, O_p^n(\nu \mid \mu)
    - \sum_{\mu = 1}^{r} \sum_{\nu = 0}^{\mu} 
    (\mu - \nu)(p - \mu)\, O_p^n(\nu \mid \mu).
\end{align*}
Simplifying the above gives
\begin{align*}
    |E_2(H_p^n)| &= 
    \frac{1}{2} \sum_{\mu = 1}^{r} \sum_{\nu = 0}^{\mu} 
    O_p^n(\nu \mid \mu)
    \left( 
        \binom{p}{2} - \frac{(p - \mu)(p - \mu - 1)}{2} 
        - \nu(p - \mu) - 2(\mu - \nu)(p - \mu)
    \right) \\
    &= 
    \frac{1}{2} \binom{p}{2} 
    \sum_{\mu = 1}^{r} \sum_{\nu = 0}^{\mu} 
    O_p^n(\nu \mid \mu)
    - 
    \frac{1}{4} \sum_{\mu = 1}^{r} (p - \mu) 
    \sum_{\nu = 0}^{\mu} 
    (p + 3\mu - 2\nu - 1)\, O_p^n(\nu \mid \mu).
\end{align*}

Finally, applying Lemmas~\ref{lem:sum-O=pn} and  ~\ref{cor:split-Opnu}, we obtain
\[
|E_2(H_p^n)| 
= \frac{1}{2} \binom{p}{2} p^n 
- \frac{1}{4} \sum_{\mu = 1}^{r} (p - \mu) 
\sum_{\nu = 0}^{\mu} 
(p + 3\mu - 2\nu - 1)\, O_p^n(\nu \mid \mu),
\]
as claimed.
\end{proof}

\noindent
From these expressions, several corollaries and degree-based results follow naturally.  
We now summarize them for completeness.

\begin{corollary}\label{coro:another_espression_of_E2}  
An equivalent form of Equation~\eqref{eq:E2} is given by
\begin{equation}
    \begin{aligned}
|E_2(H_p^n)|&=\dfrac{1}{4} \sum_{\mu = 1}^{r}  \sum_{\nu = 0}^{\mu} \left( 2p(\nu- \mu)  + \mu( 3\mu- 2\nu - 1) \right) O_p^n(\nu \mid \mu)\\
&=\dfrac{1}{4} \sum_{\mu = 1}^{r}  \sum_{\nu = 0}^{\mu} \left( 2p\nu- 2p\mu  + 3\mu^2- 2\mu\nu - \mu \right) O_p^n(\nu \mid \mu).
\end{aligned}
\end{equation}
Moreover,
\begin{equation}
    \begin{aligned}
|E_1(H_p^n)|+ |E_2(H_p^n)|
&= \dfrac{1}{4} \sum_{\mu = 1}^{r}  \sum_{\nu = 0}^{\mu} \left( 3\mu^{2} -2p\mu-\mu+4p\nu-4\mu\nu \right) O_p^n(\nu \mid \mu)\\
&= \dfrac{1}{4} \sum_{\mu = 1}^{r}  \sum_{\nu = 0}^{\mu} \left( 4\nu(p-\mu)+\mu(3\mu-2p-1)  \right) O_p^n(\nu \mid \mu).
\end{aligned}
\end{equation}
\end{corollary}

\noindent
Finally, we connect this classification to the degree structure of $H_p^n$.  
The following propositions enumerate edges according to the degrees of their endpoints, distinguishing equal-degree and off-diagonal cases.

\subsubsection*{Degree-based decomposition of the edge set}

Having classified the edges of $H_p^n$ according to the type of move, we now turn to a finer analysis based on the \emph{degrees} of their endpoints.  We denote by $\Delta(H_p^n)=f_{p,n}(p)=\binom{p}{2}$ the maximum vertex degree of $H_p^n$.

\begin{proposition}\label{prop:eqdeg-lower}
For all $p \in \mathbb{N}_3$ and $n \in \mathbb{N}_0$,  let $r = \min\{n, p\}$,  $\lambda < \Delta(H_p^n)$, and let $\mu=f^{-1}(\lambda)$.
The number of edges $\{u,v\}\in E(H_p^n)$ such that
\[
\deg(u)=\deg(v)=\lambda
\]
is given by
\begin{equation}
    \sigma(\lambda)
    = \dfrac{1}{4}   \sum_{\nu = 0}^{\mu} 
       \bigl( 3\mu^{2} - 2p\mu - \mu + 4p\nu - 4\mu\nu \bigr)
       O_p^n(\nu \mid \mu).
\end{equation}
\end{proposition}

\begin{proof}
Since $\lambda < \Delta(H_p^n)$, the inverse image $\mu=f^{-1}(\lambda)$ is unique.  
From Corollary~\ref{coro:another_espression_of_E2}, we know that the sum of $|E_1(H_p^n)|$ and $|E_2(H_p^n)|$ accounts precisely for all edges joining vertices of equal degree (i.e., those for which $o(u)=o(v)=\mu$).  
Substituting this relation gives the stated expression for $\sigma(\lambda)$.
\end{proof}

\noindent
The next proposition characterizes the edges whose endpoints both have the \emph{maximum degree}, i.e., $\lambda=\Delta(H_p^n)$.

\begin{proposition}\label{prop:eqdeg-top}  
For all $p \in \mathbb{N}_3$ and $n \in \mathbb{N}_0$,  let $r = \min\{n, p\}$ and  $\lambda=\Delta(H_p^n)=\binom{p}{2}$.
The number of edges $\{u,v\}\in E(H_p^n)$ satisfying $\deg(u)=\deg(v)=\lambda$ is
\begin{equation}\label{eq:sigma-top}
         \sigma(\lambda)=\begin{cases}
             \dfrac{1}{4}  \bigl(p^{2} -p\bigr)\, O_p^n(p), & \text{if } o(u)=o(v)=p,\\[4pt]
             \dfrac{1}{4}   \displaystyle\sum_{\nu = 0}^{p-1} 
             \bigl(p^{2}-5p+4\nu+4\bigr) O_p^n(\nu \mid p-1),
             & \text{if } o(u)=o(v)=p-1,\\[8pt]
             \displaystyle\sum_{\nu = 0}^{p-1}  (p-\nu-1)\, O_p^n(\nu \mid p-1),
             & \text{if } o(u)=p-1,\; o(v)=p.
         \end{cases}
\end{equation}
Consequently, we obtain the compact identity
\begin{equation}
        \sigma(\lambda)
        = \dfrac{1}{4}\,p
          \begin{Bmatrix}
             n+1\\p
          \end{Bmatrix}.
\end{equation}
\end{proposition}

\begin{proof}
If $\deg(u)=\deg(v)=\lambda=\Delta(H_p^n)$, then both endpoints have either all pegs occupied or all but one occupied. Hence,
\[
\{o(u),o(v)\}\in\bigl\{\{p,p\},\{p-1,p-1\},\{p-1,p\}\bigr\}.
\]
We analyze these cases separately.

\noindent
\begin{itemize}
    \item Case 1, $o(u)=o(v)=p$:

Here $\mu=f^{-1}(\lambda)=p$.  
Using Corollary~\ref{coro:another_espression_of_E2} with $\mu=p$, we find
\begin{align*}
\sigma(\lambda)
&= \dfrac{1}{4}   \sum_{\nu = 0}^{r} 
   \bigl(3p^{2} -2pp -p + 4p\nu -4p\nu\bigr) O_p^n(\nu \mid p)\\
&= \dfrac{1}{4}  (p^{2}-p)\,O_p^n(p).
\end{align*}

\item Case 2, $o(u)=o(v)=p-1$: 

Now $\mu=p-1$. Substituting $\mu=p-1$ in the same expression gives
\begin{align*}
\sigma(\lambda)
&= \dfrac{1}{4}   \sum_{\nu = 0}^{p-1} 
   \bigl(3(p-1)^2 -2p(p-1)-(p-1)+4p\nu-4\nu(p-1)\bigr)
   O_p^n(\nu \mid p-1)\\
&= \dfrac{1}{4}   \sum_{\nu = 0}^{p-1} 
   \bigl(p^{2}-5p+4\nu+4\bigr)\,O_p^n(\nu \mid p-1).
\end{align*}

\item  Case 3,  $\{o(u),o(v)\}=\{p-1,p\}$:

Here we use the expression from case~\texttt{(c\textsubscript{3})}, corresponding to a move from a multiton to an empty peg:
\begin{align*}
\sigma(\lambda)
&= \sum_{\nu = 0}^{p-1}(p-\mu)(\mu-\nu)\,O_p^n(\nu \mid \mu)
 = \sum_{\nu = 0}^{p-1}(p-(p-1))((p-1)-\nu)\,O_p^n(\nu \mid p-1)\\
&= \sum_{\nu = 0}^{p-1}(p-\nu-1)\,O_p^n(\nu \mid p-1).
\end{align*}
\end{itemize}
 
Summing the three contributions and applying the Stirling-type identity
\[
p\begin{Bmatrix} n\\p \end{Bmatrix}
   + \begin{Bmatrix} n\\p-1 \end{Bmatrix}
   = \begin{Bmatrix} n+1\\p \end{Bmatrix},
\]
we obtain the compact form

\begingroup
\allowdisplaybreaks
     \begin{align*}
         \sigma(\lambda)&= \dfrac{1}{4}  \left( p^{2} -p  \right) O_p^n(  p)+ \dfrac{1}{4}   \sum_{\nu = 0}^{p-1} \left( p^{2}-5p+4\nu+4 \right) O_p^n(\nu \mid p-1)+\sum_{\nu = 0}^{p-1}  (p-\nu-1) O_p^n(\nu \mid p-1)\\
         &=\dfrac{1}{4}  \left( p^{2} -p  \right) O_p^n(  p)+\dfrac{1}{4}  \left( p^{2} -p  \right)  \sum_{\nu = 0}^{p-1}   O_p^n(\nu \mid p-1)\\
         &=\dfrac{1}{4}  \left( p^{2} -p  \right) O_p^n(  p)+\dfrac{1}{4}  \left( p^{2} -p  \right)    O_p^n(  p-1)\\
         &=\dfrac{1}{4}  \left( p^{2} -p  \right) (O_p^n(  p)+ O_p^n(  p-1))\\
         &=\dfrac{1}{4}  \left( p^{2} -p  \right) \left( \begin{Bmatrix}
             n\\p
         \end{Bmatrix}p! +\begin{Bmatrix}
             n\\p-1
         \end{Bmatrix}(p-1)! \right)\\
         &=\dfrac{1}{4}  \left( p -1  \right) p! \left( p\begin{Bmatrix}
             n\\p
         \end{Bmatrix} +\begin{Bmatrix}
             n\\p-1
         \end{Bmatrix} \right)\\
         &=\dfrac{1}{4}  p\begin{Bmatrix}
             n+1\\p
         \end{Bmatrix}.\qedhere
     \end{align*}
\endgroup
\end{proof}

\noindent
We next consider \emph{off-diagonal edges}, i.e., edges connecting vertices of distinct degrees. These correspond precisely to edges where the number of occupied pegs differs by one between the two endpoints.

\begin{proposition}\label{prop:offdiag-general}  
For all $p \in \mathbb{N}_3$ and $n \in \mathbb{N}_0$,  let   $\lambda < \Delta(H_p^n)$ and  $\mu=f^{-1}(\lambda)$.  
The number of edges $\{u,v\}\in E(H_p^n)$ such that $\deg(u)<\deg(v)=\lambda$ is
\begin{equation}
    \sigma(\lambda)
    = \sum_{\nu = 0}^{\mu} (p-\mu)(\mu-\nu)\,O_p^n(\nu \mid \mu).
\end{equation}
\end{proposition}

\begin{proof}
If $\deg(u)<\deg(v)=\lambda$, then $o(u)<o(v)$ and necessarily $o(v)=o(u)+1=\mu+1$.  
Since $\mu\le p-2$, this correspondence between $\mu$ and $\lambda$ is one-to-one.  
The number of such edges equals the number of moves of type \texttt{(c\textsubscript{3})}, i.e., moves from a multiton peg to an empty peg in a state with $\mu$ occupied pegs.  
Each such move can be made by choosing one of the $(\mu-\nu)$ multitons and one of the $(p-\mu)$ empty pegs, giving the stated formula.
\end{proof}

\noindent
When the upper endpoint reaches the maximal degree, two additional contributions appear, corresponding to $o(v)=p-1$ and $o(v)=p$.

\begin{proposition}\label{prop:offdiag-top}
For all $p \in \mathbb{N}_3$ and $n \in \mathbb{N}_0$,  let   $\lambda=\Delta(H_p^n)=\binom{p}{2}$.  
The number of edges $\{u,v\}\in E(H_p^n)$ such that $\deg(u)<\deg(v)=\lambda$ is
\begin{equation}
    \sigma(\lambda)
    = \sum_{\nu = 0}^{p-1} (p-1-\nu)\,O_p^n(\nu \mid p-1)
      + \sum_{\nu = 0}^{p-2} 2(p-2-\nu)\,O_p^n(\nu \mid p-2).
\end{equation}
\end{proposition}

\begin{proof}
If $\deg(v)=\lambda=\Delta(H_p^n)$, then $o(v)\in\{p-1,p\}$.
\begin{itemize}
    \item Case 1, $o(v)=p$:   
Then $o(u)=p-1$, and by applying case~\texttt{(c\textsubscript{3})} with $\mu=p-1$, we have
\[
\sigma(\lambda)
= \sum_{\nu = 0}^{p-1}(p-(p-1))((p-1)-\nu)\,O_p^n(\nu \mid p-1)
= \sum_{\nu = 0}^{p-1}(p-1-\nu)\,O_p^n(\nu \mid p-1).
\]
    \item Case 2, $o(v)=p-1$:  
Then $o(u)=p-2$, and again by case~\texttt{(c\textsubscript{3})}, we have
\[
\sigma(\lambda)
= \sum_{\nu = 0}^{p-2}(p-(p-2))((p-2)-\nu)\,O_p^n(\nu \mid p-2)
= \sum_{\nu = 0}^{p-2}2(p-2-\nu)\,O_p^n(\nu \mid p-2).
\]
\end{itemize}

Adding both contributions completes the proof.
\end{proof}

\noindent
Gathering all the above results, we can now state a unified expression that encompasses both the general and top-degree cases.

\begin{proposition}\label{prop:offdiag-lambda}
For all $p \in \mathbb{N}_3$ and $n \in \mathbb{N}_0$. For every $\lambda\le\Delta(H_p^n)$, let $\mu=f^{-1}(\lambda)$.  
The number of edges $\{u,v\}\in E(H_p^n)$ satisfying $\deg(u)<\deg(v)=\lambda$ is
\begin{equation}
    \sigma(\lambda)
    = \sum_{\nu = 0}^{\mu} (p-\mu)(\mu-\nu)\,O_p^n(\nu \mid \mu).
\end{equation}
\end{proposition}

\begin{proof}
The formula of Proposition~\ref{prop:offdiag-general} already holds for all $\lambda<\Delta(H_p^n)$.  
For $\lambda=\Delta(H_p^n)$, the two contributions given in Proposition~\ref{prop:offdiag-top} are consistent with the same expression, since for $\mu=p-1$ the term $(p-\mu)=1$, and for $\mu=p-2$ we obtain the multiplicative factor $2$ accounting for the two admissible target configurations.  
Thus, the formula extends uniformly to all $\lambda\le\Delta(H_p^n)$.
\end{proof}

\noindent
This completes the degree-based decomposition of the edge set of $H_p^n$, providing an explicit enumeration of all edge classes as a function of the occupancy parameters $(\mu,\nu)$ and the associated degree values.

\section{Explicit expression of the M-polynomial }
\label{sec:explicit}

Let $r=\min\{n,p\}$ and define the (multi)set of realized degree values
\[
\Gamma_p^n \;=\; \Big\{ \lambda_\mu : \mu\in[r] \Big\}
\qquad\text{with}\qquad
\lambda_\mu \;=\; f_{p,n}(\mu)\,.
\]
By Proposition~\ref{prop:deg-noninj} the map $\mu\mapsto\lambda_\mu$ is bijective on $[r]\setminus\{p\}$
and satisfies $\lambda_{p-1}=\lambda_p=\Delta(H_p^n)$ when $r=p$ (top degree collision).
For notational convenience we set
\[
\Delta = \Delta(H_p^n)=\lambda_{p-1}=\lambda_p
\quad\text{(when $r=p$)},\qquad
\delta = \delta(H_p^n)=\lambda_1=p-1.
\]

The M-polynomial of $H_p^n$ is
\[
M(H_p^n;x,y)
\;=\;
\sum_{i\le j} m_{i,j}(H_p^n)\,x^i y^j
\;=\;
\sum_{\lambda\le \lambda'} m_{\lambda,\lambda'}(H_p^n)\,x^\lambda y^{\lambda'}\,,
\]
where $m_{\lambda,\lambda'}(H_p^n)$ counts edges $\{u,v\}\in E(H_p^n)$ with
$\{\deg(u),\deg(v)\}=\{\lambda,\lambda'\}$.

Our strategy is to \emph{index} degree classes by the occupancy parameter $\mu$, and to use the
edge partition of Theorem~\ref{thm:edge_partition}  together with the occupancy
enumerators $O_n^p(\nu\mid \mu)$ (Proposition~\ref{th:refined}).  This yields closed forms
for all diagonal coefficients $m_{\lambda,\lambda}$ and all off–diagonal coefficients
between \emph{adjacent} degree classes (which are the only nonzero off–diagonal terms by
Proposition~\ref{prop:deg-diff}).

\subsection*{Degree classes and adjacency pattern}

For $1\le \mu\le r$, define the \emph{$\mu$–degree class}
\[
\mathcal{D}_\mu \;=\;\{\, s\in V(H_p^n) \;:\; o(s)=\mu \,\},
\qquad
\deg(\mathcal{D}_\mu)=\lambda_\mu.
\]
By Proposition~\ref{prop:deg-diff} and Corollary~\ref{cor:deg-step},
if $\{u,v\}\in E(H_p^n)$ with $o(u)\le o(v)$, then $|o(u)-o(v)|\in\{0,1\}$, hence
\begin{equation}
    \{\deg(u),\deg(v)\}\in
\begin{cases}
\{\lambda_\mu\} & \text{if } o(u)=o(v)=\mu,\\[2mm]
\{\lambda_{\mu-1},\lambda_\mu\} & \text{if } o(u)=\mu-1,\;o(v)=\mu.
\end{cases}
\end{equation}

Therefore the only possibly nonzero off–diagonal coefficients are those between \emph{consecutive}
occupancy classes.

\subsection*{Diagonal coefficients $m_{\lambda,\lambda}$}

\begin{lemma}[Equal-degree edges inside a fixed occupancy]
\label{lem:eqdeg-lower}
For all $p \in \mathbb{N}_3$ and $n \in \mathbb{N}_0$,  let $r = \min\{n, p\}$ and  $1\le \mu\le \min\{r,p-2\}$. Then
\begin{equation}
    m_{\lambda_\mu,\lambda_\mu}(H_p^n)
\;=\;
\frac14\sum_{\nu=0}^{\mu}
\Big(3\mu^2-2p\mu-\mu+4p\nu-4\mu\nu\Big)\; O_n^p(\nu\mid \mu).
\end{equation}
\end{lemma}

\begin{proof}
By Proposition~\ref{prop:eqdeg-lower}, for any degree value
$\lambda<\Delta$ there is a unique $\mu=f_{p,n}^{-1}(\lambda)$, and the number of equal–degree
edges with both endpoints in $\mathcal{D}_\mu$ is exactly
\[
\frac14\sum_{\nu=0}^{\mu}
\Big(3\mu^2-2p\mu-\mu+4p\nu-4\mu\nu\Big)\; O_n^p(\nu\mid \mu).
\]
Since here $\lambda_\mu<\Delta$ for $\mu\le p-2$, the claim follows.
\end{proof}

\begin{lemma}[Equal-degree edges at top degree]
\label{lem:eqdeg-top}
For all $p \in \mathbb{N}_3$ and $n \in \mathbb{N}_0$. Assume $r=p$ so that $\lambda_{p-1}=\lambda_p=\Delta$.
Then
\begin{equation}
    m_{\Delta,\Delta}(H_p^n)
\;=\;
\frac14\Big(p^2-p\Big)\Big( O_n^p(p)+O_n^p(p-1)\Big).
\end{equation}
Equivalently, using $O_n^p(\mu)=\begin{Bmatrix} n \\ \mu \end{Bmatrix}\,p^{\underline{\mu}}$,
\begin{equation}
    m_{\Delta,\Delta}(H_p^n)
\;=\;
\frac14p!\Big(p^2-p\Big)\Big(
\begin{Bmatrix} n \\ p \end{Bmatrix}
\;+\;
\begin{Bmatrix} n \\ p-1 \end{Bmatrix}
\Big).
\end{equation}
\end{lemma}

\begin{proof}
By Proposition~\ref{prop:eqdeg-top}, edges with both endpoints of
degree $\Delta$ come from three disjoint configurations:
(i) $o(u)=o(v)=p$;
(ii) $o(u)=o(v)=p-1$; and
(iii) $\{o(u),o(v)\}=\{p-1,p\}$.
Adding the three expressions given in Proposition~\ref{prop:eqdeg-top} and simplifying using
Corollary~\ref{cor:split-Opnu} yields
\[
m_{\Delta,\Delta}
= \frac14\,(p^2-p)O_n^p(p) \;+\; \frac14\,(p^2-p)O_n^p(p-1)
= \frac14\,(p^2-p)\Big(O_n^p(p)+O_n^p(p-1)\Big),
\]
as claimed.
\end{proof}

\subsection*{Off–diagonal coefficients $m_{\lambda,\lambda'}$ between adjacent classes}

\begin{lemma}[Edges between consecutive occupancies]
\label{lem:offdiag-general}
For all $p \in \mathbb{N}_3$ and $n \in \mathbb{N}_0$,  let $r = \min\{n, p\}$ and  $2\le \mu\le r$. Then the number of edges with
$o(u)=\mu-1$, $o(v)=\mu$ (hence $\{\deg(u),\deg(v)\}=\{\lambda_{\mu-1},\lambda_\mu\}$)
equals
\begin{equation}
    m_{\lambda_{\mu-1},\lambda_\mu}(H_p^n)
\;=\;
\sum_{\nu=0}^{\mu-1}\big(p-\mu+1\big)\,(\mu-1-\nu)\; O_n^p(\nu\mid \mu-1).
\end{equation}
\end{lemma}

\begin{proof}
This is exactly the case $(\texttt{c}_{3})$) of Theorem~\ref{thm:edge_partition} with $\mu$ replaced by $\mu-1$,
counting moves from a multiton peg to an empty peg.  Proposition~\ref{prop:offdiag-lambda}  gives for $\lambda=\lambda_\mu$:
\[
\sigma(\lambda_\mu)
=
\sum_{\nu=0}^{\mu-1} (p-(\mu-1))\,\big((\mu-1)-\nu\big)\; O_n^p(\nu\mid \mu-1),
\]
which is precisely $m_{\lambda_{\mu-1},\lambda_\mu}$.
\end{proof}

\begin{lemma}[Off–diagonal edges to top degree]
\label{lem:offdiag-top}
For all $p \in \mathbb{N}_3$ and $n \in \mathbb{N}_0$. Assume $r=p$. Then the total number of edges with the higher endpoint at degree $\Delta$
(i.e., with $\{\deg(u),\deg(v)\}=\{\lambda,\Delta\}$ and $\lambda<\Delta$) equals
\begin{equation}
    m_{\lambda_{p-1},\Delta}(H_p^n)
\;=\;
\sum_{\nu=0}^{p-1}(p-1-\nu)\,O_n^p(\nu\mid p-1)
\;+\;
\sum_{\nu=0}^{p-2} 2\,(p-2-\nu)\,O_n^p(\nu\mid p-2).
\end{equation}
\end{lemma}

\begin{proof}
By Proposition~\ref{prop:offdiag-top}, off–diagonal edges with higher degree
$\Delta$ arise in two disjoint situations:
(i) $o(v)=p$, $o(u)=p-1$ (coefficient $\sum_{\nu=0}^{p-1}(p-1-\nu)O_n^p(\nu\mid p-1)$);
(ii) $o(v)=p-1$, $o(u)=p-2$ (coefficient $\sum_{\nu=0}^{p-2}2(p-2-\nu)O_n^p(\nu\mid p-2)$).
Adding the two contributions yields the stated formula.
\end{proof}

\subsection*{The main formula}

\begin{theorem}[Explicit M-polynomial of $H_p^n$]
\label{thm:Mpoly} For all $p \in \mathbb{N}$ and $n \in \mathbb{N}_0$. 
Let $r=\min\{n,p\}$ and $\lambda_\mu=f_{p,n}(\mu)$ for $\mu\in[r]$.
Then
\begin{equation}
    \begin{aligned}
M(H_p^n;x,y)
&=
\underbrace{\sum_{\mu=1}^{\min\{r,p-2\}}
m_{\lambda_\mu,\lambda_\mu}(H_p^n)\,x^{\lambda_\mu}y^{\lambda_\mu}}_{\text{diagonal terms below top degree}}
\;+\;
\mathbf{1}_{\{r=p\}}\; m_{\Delta,\Delta}(H_p^n)\,x^\Delta y^\Delta
\\[1mm]
&\quad
+\;\sum_{\mu=2}^{\min\{r,p-1\}}
m_{\lambda_{\mu-1},\lambda_\mu}(H_p^n)\,x^{\lambda_{\mu-1}}y^{\lambda_\mu}+\;\mathbf{1}_{\{r=p\}}\; m_{\lambda_{p-1},\Delta}(H_p^n)\,x^{\lambda_{p-1}}y^{\Delta}.
\end{aligned}
\end{equation}

where the coefficients are given explicitly by
\[
\begin{aligned}
    m_{\lambda_\mu,\lambda_\mu}(H_p^n)
=&
\frac14\sum_{\nu=0}^{\mu}
\Big(3\mu^2-2p\mu-\mu+4p\nu-4\mu\nu\Big)\; O_n^p(\nu\mid \mu)
&& (1\le \mu\le \min\{r,p-2\}),\\
m_{\Delta,\Delta}(H_p^n)
=&
\frac14\Big(p^2-p\Big)\Big( O_n^p(p)+O_n^p(p-1)\Big)
&&\text{(when $r=p$)},\\
m_{\lambda_{\mu-1},\lambda_\mu}(H_p^n)
=&
\sum_{\nu=0}^{\mu-1}\big(p-\mu+1\big)\,(\mu-1-\nu)\; O_n^p(\nu\mid \mu-1)
&& (2\le \mu\le \min\{r,p-1\}),\\
m_{\lambda_{p-1},\Delta}(H_p^n)
=&
\sum_{\nu=0}^{p-1}(p-1-\nu)\,O_n^p(\nu\mid p-1)
\; &&  \\
&+\;
\sum_{\nu=0}^{p-2} 2\,(p-2-\nu)\,O_n^p(\nu\mid p-2)
&& \text{(when $r=p$)}. 
\end{aligned}
\]
\end{theorem}

\begin{proof}
By the discussion preceding Lemma~\ref{lem:eqdeg-lower}, all diagonal coefficients below $\Delta$
are equal–degree edges within a single occupancy class $\mathcal{D}_\mu$ ($\mu\le p-2$); these are
counted by Proposition~\ref{prop:eqdeg-lower}, giving the first line of coefficients.

At top degree (when $r=p$), equal–degree edges arise from three disjoint subcases
$\{p,p\}$, $\{p-1,p-1\}$, and $\{p-1,p\}$; Proposition~\ref{prop:eqdeg-top} sums these to the stated
closed form for $m_{\Delta,\Delta}$.

For off–diagonal terms between adjacent degree classes $\lambda_{\mu-1}$ and $\lambda_\mu$ with
$2\le \mu\le p-1$, Proposition~\ref{prop:offdiag-general} (case $(\texttt{c}_{3})$) yields the stated coefficient.
Finally, when $r=p$, edges with higher degree $\Delta$ have two sources (from $p-1$ to $p$ and
from $p-2$ to $p-1$); adding these contributions as in Proposition~\ref{prop:offdiag-top} gives the
last line.

No other off–diagonal pairs occur: by Proposition~\ref{prop:deg-diff}, if $\{u,v\}$ is an edge with
$o(u)\le o(v)$, then $o(v)-o(u)\in\{0,1\}$, hence degree differences can only connect consecutive
occupancy classes.  Collecting all contributions produces the asserted decomposition.
\end{proof}

\begin{proposition}[Edge count recovered]
Summing all coefficients of $M(H_p^n;x,y)$ at $x=y=1$ yields
\[
\sum_{\lambda\le \lambda'} m_{\lambda,\lambda'}(H_p^n) \;=\; |E(H_p^n)|,
\]
and the right–hand side simplifies to the closed forms of
Theorem~\ref{theo:edges-sum-forms} (equations~\eqref{eq:Ehpn-sum-mu} and~\eqref{eq:Ehpn-alt}).
\end{proposition}

\begin{proof}
Directly from the definition of $m_{\lambda,\lambda'}$ and the partition of the edge set by
Propositions~\ref{prop:two-block} and~\ref{thm:edge_partition}; evaluating $M$ at $(1,1)$ counts each
edge exactly once.
\end{proof}

\begin{remark}[Closed forms in Stirling numbers]
Using Proposition~\ref{th:Opmn} we may rewrite $O_n^p(\nu\mid \mu)$ as
\[
O_n^p(\nu\mid \mu)
\;=\;
\binom{n}{\nu}\,
\begin{Bmatrix} n-\nu \\ \mu-\nu \end{Bmatrix}_{\ge 2}\, p^{\underline{\mu}},
\qquad
p^{\underline{\mu}} = p(p-1)\cdots (p-\mu+1),
\]
so Theorem~\ref{thm:Mpoly} becomes a fully explicit expression in terms of (associated) Stirling
numbers and falling factorials.  In particular,
\[
O_n^p(\mu)=\sum_{\nu=0}^{\mu} O_n^p(\nu\mid \mu)
=
\begin{Bmatrix} n \\ \mu \end{Bmatrix}\,p^{\underline{\mu}},
\]
and the identity
\(
\sum_{\mu=1}^{r} O_n^p(\mu) = p^n
\)
recovers Lemma~\ref{lem:sum-O=pn}.
\end{remark}

\begin{remark}[Adjacency constraint]
     By Proposition~\ref{prop:deg-diff} and Corollary~\ref{cor:deg-step},
edges only connect equal occupancies ($\mu\to\mu$) or \emph{adjacent} occupancies ($\mu\to\mu+1$).
Therefore, in $M(H_p^n;x,y)$ only diagonal terms $x^{\lambda_\mu}y^{\lambda_\mu}$ and consecutive off–diagonal terms
$x^{\lambda_\mu}y^{\lambda_{\mu+1}}$ appear (with $i\le j$), exactly as enforced in Theorem~\ref{thm:Mpoly}.
\end{remark}

\subsection*{Worked examples and index extraction}

We illustrate Theorem~\ref{thm:Mpoly} with two concrete instances: the classical three–peg,
three–disc graph $H_3^3$, and the four–peg, two–disc graph $H_4^2$. In each case we (i) compute
the occupancy enumerators $O_n^p(\nu\mid\mu)$, (ii) assemble the edge counts by the
$(\texttt{c}_{1})$--$(\texttt{c}_{3})$ decomposition of Theorem~\ref{thm:edge_partition}, and (iii) synthesize
$M(H_p^n;x,y)$. Finally we verify by extracting $M_1$ and $M_2$ via the operator table.

\paragraph{Example 1: $H_3^3$.}
Here $p=3$, $n=3$, $r=\min\{n,p\}=3$; the realized degrees are
\[
\lambda_1=f_{3,3}(1)=2,\qquad
\lambda_2=f_{3,3}(2)=3,\qquad
\lambda_3=f_{3,3}(3)=3=\Delta,
\]
so $\Gamma_3^3=\{2,3\}$ and the top degree is $\Delta=3$ (with a collision at $\mu=2,3$).

\emph{Occupancy counts.}
Using Proposition~\ref{th:refined},
\[
O_3^3(\nu\mid\mu)=\binom{3}{\nu}\,
\begin{Bmatrix} 3-\nu \\ \mu-\nu \end{Bmatrix}_{\ge 2}\,
(3)_\mu,
\qquad (3)_1=3,\ (3)_2=6,\ (3)_3=6.
\]
A direct check gives
\[
\begin{array}{c|ccc}
\mu\backslash \nu & 0 & 1 & 2,3\\\hline
1 & 3 & 0 & 0\\
2 & 0 & 18 & 0\\
3 & 0 & 0 & 6
\end{array}
\qquad\Longrightarrow\qquad
O_3^3(1)=3,\; O_3^3(2)=18,\; O_3^3(3)=6,
\]
and $O_3^3(1)+O_3^3(2)+O_3^3(3)=27=3^3$.

\emph{Edges between occupancy classes.}
By $(\texttt{c}_{3})$, edges from $\mu$ to $\mu+1$ are counted (with orientation) by
$(\mu-\nu)(p-\mu)$ from a state of type $(\nu\mid\mu)$; for undirected edges and fixed
$\mu$ this is a \emph{single} count (no symmetry division because $\mu\neq \mu+1$).

\smallskip
\noindent
$\bullet$ Between $\mu=1$ (all three discs on one peg: $\nu=0$) and $\mu=2$:
each such state contributes $(\mu-\nu)(p-\mu)=1\cdot 2=2$ edges; there are $O_3^3(0\mid 1)=3$
such states, hence
\[
|E(\mu=1\leftrightarrow \mu=2)| \;=\; 3\cdot 2 \;=\; 6.
\]
These are edges of \emph{degree pair} $(2,3)$, thus they contribute to $m_{2,3}$.

\smallskip
\noindent
$\bullet$ Between $\mu=2$ (necessarily $\nu=1$) and $\mu=3$ (necessarily $\nu=3$):
each $\mu=2$ state contributes $(\mu-\nu)(p-\mu)=1\cdot 1=1$; there are $O_3^3(1\mid 2)=18$
such states, hence
\[
|E(\mu=2\leftrightarrow \mu=3)| \;=\; 18.
\]
Here \emph{both} endpoints have degree $3$ (since $\lambda_2=\lambda_3=3$), so these edges
contribute to $m_{3,3}$.

\emph{Edges within occupancy classes.}
For $\mu=1$ there are no within-class edges (only one occupied peg, so any legal move
increases $\mu$). For $\mu=3$ there are again no within-class edges:
from a state with three singletons, any legal move decreases the number of occupied pegs.
Thus the only within-class contribution is for $\mu=2$.

At $\mu=2$ we necessarily have $\nu=1$ (one singleton, one multiton). Edges that keep $\mu=2$
are of types $(\texttt{c}_{1})$ SP$\to$EP and $(\texttt{c}_{2})$ MP$\to$OP. The $(\texttt{c}_{1})$ contribution depends only on
$(\nu\mid\mu)$, but $(\texttt{c}_{2})$ depends on the top-disc size ordering at both pegs and is not
directly multiplicative in $(\nu,\mu)$.
Instead of a delicate direct count, we now \emph{balance} using the total edge count. Indeed,
\[
|E(H_3^3)| \;=\; \frac12\sum_{s\in V}\deg(s)
= \frac12\big( O_3^3(1)\cdot 2 + (O_3^3(2)+O_3^3(3))\cdot 3\big)
=\frac12(3\cdot 2 + 24\cdot 3)=39.
\]
Subtracting the already-counted cross–class edges (6 between $\mu=1,2$ and 18 between
$\mu=2,3$) leaves
\[
|E(\text{within }\mu=2)| \;=\; 39-(6+18)=15,
\]
and these are all of degree pair $(3,3)$.

\emph{M-polynomial.}
Collecting all contributions,
\[
m_{2,2}=0,\qquad m_{2,3}=6,\qquad m_{3,3}=18+15=33,
\]
and therefore
\[
M(H_3^3;x,y)\;=\;6\,x^2y^3\;+\;33\,x^3y^3.
\]

\emph{Checks via indices.}
Applying the operators at $(x,y)=(1,1)$:
\[
M_1(H_3^3)=(D_x+D_y)M\big|_{1,1}
=6(2+3)+33(3+3)=30+198=228.
\]
This equals $\sum_{v\in V}\deg(v)^2=3\cdot 2^2+24\cdot 3^2=12+216=228$.
Also
\[
M_2(H_3^3)=(D_xD_y)M\big|_{1,1}=6(2\cdot 3)+33(3\cdot 3)=36+297=333.
\]
Finally, $M(1,1)=6+33=39=|E(H_3^3)|$.

\paragraph{Example 2: $H_4^2$.}
Here $p=4$, $n=2$, $r=2$; realized degrees are
\[
\lambda_1=f_{4,2}(1)=3,\qquad \lambda_2=f_{4,2}(2)=5,
\]
so $\Gamma_4^2=\{3,5\}$.

\emph{Occupancy counts.}
Using $O_2^4(\nu\mid\mu)=\binom{2}{\nu}\,
\begin{Bmatrix} 2-\nu \\ \mu-\nu \end{Bmatrix}_{\ge 2}\,(4)_\mu$ with
$(4)_1=4$, $(4)_2=12$, we obtain
\[
O_2^4(0\mid 1)=4,\quad O_2^4(1\mid 1)=0,\qquad
O_2^4(2\mid 2)=12,\quad O_2^4(0\mid 2)=O_2^4(1\mid 2)=0,
\]
hence $O_2^4(1)=4$, $O_2^4(2)=12$, and $4+12=16=4^2$.

\emph{Edges between classes.}
By $(\texttt{c}_{3})$, from $\mu=1$ ($\nu=0$) to $\mu=2$:
each state contributes $(\mu-\nu)(p-\mu)=1\cdot 3=3$ moves; there are $4$ such states, hence
\[
|E(\mu=1\leftrightarrow \mu=2)| \;=\; 4\cdot 3 \;=\; 12.
\]
These have degree pair $(3,5)$ and contribute to $m_{3,5}$.

\emph{Edges within classes.}
Within $\mu=1$ there are none (any move increases occupancy).
Within $\mu=2$ we are in the case $\nu=2$ (two singletons). The only moves that keep $\mu=2$
are of type $(\texttt{c}_{1})$ SP$\to$EP (since $(\texttt{c}_{2})$ needs a multiton, which is absent).
Each $\nu=2$ state has $\nu(p-\mu)=2\cdot 2=4$ such moves; dividing by $2$ to remove symmetry
yields $2$ undirected edges per state, hence
\[
|E(\mu=2 \text{ within})| \;=\; 12\cdot 2 \;=\; 24,
\]
all with degree pair $(5,5)$.

\emph{Total and M-polynomial.}
The total number of edges is
\[
|E(H_4^2)|=\frac12\big( O_2^4(1)\cdot 3 + O_2^4(2)\cdot 5\big)
=\frac12(4\cdot 3+12\cdot 5)=36,
\]
in agreement with $12+24=36$ counted above.
Therefore
\[
M(H_4^2;x,y)\;=\;12\,x^3y^5\;+\;24\,x^5y^5.
\]

\emph{Checks via indices.}
\[
M_1(H_4^2)=12(3+5)+24(5+5)=96+240=336
=\sum_{v}\deg(v)^2=4\cdot 3^2+12\cdot 5^2=36+300=336,
\]
\[
M_2(H_4^2)=12(3\cdot 5)+24(5\cdot 5)=180+600=780,
\qquad
M(1,1)=12+24=36=|E(H_4^2)|.
\]

\begin{remark}[How to automate larger instances]
For general $(p,n)$, Theorem~\ref{thm:Mpoly} reduces the computation of $M(H_p^n;x,y)$ to:
\begin{enumerate}
\item tabulate $O_n^p(\nu\mid \mu)$ for $1\le \mu\le r=\min\{n,p\}$ and $0\le \nu\le \mu$ using
\(
O_n^p(\nu\mid \mu)=\binom{n}{\nu}
\begin{Bmatrix} n-\nu \\ \mu-\nu \end{Bmatrix}_{\ge 2}\, p^{\underline{\mu}},
\)
\item evaluate the coefficient sums in Theorem~\ref{thm:Mpoly} for
$m_{\lambda_\mu,\lambda_\mu}$ and $m_{\lambda_{\mu-1},\lambda_\mu}$,
\item assemble $M$ and, if desired, apply the differential/integral operators of Table~\ref{tab:indices}
to extract the indices in closed form.
\end{enumerate}
This workflow is quasi–linear in $n$ for fixed $p$ assuming precomputed associated Stirling numbers. 

This scheme is illustrated through the three sequential Tables~\ref{tab:stirling}, \ref{tab:O4n-valus}, and~\ref{tab:MH4n-coeffs}. 
Table~\ref{tab:MH4n-coeffs} presents the coefficients of $M(H_4^{n}; x, y)$ for $n = 1, \dots, 8$, which are computed from the values of $O_{4}^{n}(\nu \mid \mu)$ shown in Table~\ref{tab:O4n-valus}. 
These values, in turn, are obtained using the associated Stirling numbers provided in Table~\ref{tab:stirling}.

\end{remark}

\begin{table}[H]
\centering
\renewcommand{\arraystretch}{1.2}
\begin{tabular}{@{}l|rrrrrrrrr@{}}
\toprule
\diagbox[height=0.75cm]{ $j$ }{$i$} & 0 & 1 & 2 & 3 & 4 & 5 & 6 & 7 & 8 \\ \midrule
0 & 1 & 0 & 0 & 0 & 0 & 0 & 0 & 0 & 0 \\
1 & 0 & 0 & 1 & 1 & 1 & 1 & 1 & 1 & 1 \\
2 & 0 & 0 & 0 & 0 & 3 & 10 & 25 & 52 & 105 \\
3 & 0 & 0 & 0 & 0 & 0 & 0 & 15 & 70 & 280 \\
4 & 0 & 0 & 0 & 0 & 0 & 0 & 0 & 0 & 105 \\ \bottomrule
\end{tabular}
\caption{Associated Stirling numbers of the second kind 
$\bigl\{\!\!\begin{smallmatrix} i\\ j\end{smallmatrix}\!\!\bigr\}_{\ge 2}$ 
used in the computation of $O_4^n(\nu\mid\mu)$ for $n\le 8$, with $i=n-\nu$ and $j=\mu-\nu$.}\label{tab:stirling}
\end{table}

\begin{table}[H]
\centering

\renewcommand{\arraystretch}{1.2}
\begin{tabular}{c|rrrrrrrr}
$(\mu,\nu)$ & $n=1$ & $2$ & $3$ & $4$ & $5$ & $6$ & $7$ & $8$ \\ \hline
(1,0) & 0 & 4 & 4 & 4 & 4 & 4 & 4 & 4 \\
(1,1) & 4 & 0 & 0 & 0 & 0 & 0 & 0 & 0 \\
(2,0) & 0 & 0 & 0 & 36 & 120 & 300 & 672 & 1428 \\
(2,1) & 0 & 0 & 36 & 48 & 60 & 72 & 84 & 96 \\
(2,2) & 0 & 12 & 0 & 0 & 0 & 0 & 0 & 0 \\
(3,0) & 0 & 0 & 0 & 0 & 0 & 360 & 2520 & 11760 \\
(3,1) & 0 & 0 & 0 & 0 & 360 & 1440 & 4200 & 10752 \\
(3,2) & 0 & 0 & 0 & 144 & 240 & 360 & 504 & 672 \\
(3,3) & 0 & 0 & 24 & 0 & 0 & 0 & 0 & 0 \\
(4,0) & 0 & 0 & 0 & 0 & 0 & 0 & 0 & 2520 \\
(4,1) & 0 & 0 & 0 & 0 & 0 & 0 & 2520 & 20160 \\
(4,2) & 0 & 0 & 0 & 0 & 0 & 1080 & 5040 & 16800 \\
(4,3) & 0 & 0 & 0 & 0 & 240 & 480 & 840 & 1344 \\
(4,4) & 0 & 0 & 0 & 24 & 0 & 0 & 0 & 0 \\
\end{tabular}
\caption{Values of $O_{4}^{n}(\nu\mid\mu)$ for $1\le n\le 8$.}\label{tab:O4n-valus}
\end{table}

\begin{table}[H]
\centering

\begin{tabular}{r|rrrrrr}
\toprule
$n$ & $m_{3,3}$ & $m_{3,5}$ & $m_{5,5}$ & $m_{5,6}$ & $m_{6,6}$ & $|E(H_4^n)|$ \\
\midrule
1 & 6 & 0 & 0 & 0 & 0 & 6 \\
2 & 0 & 12 & 24 & 0 & 0 & 36 \\
3 & 0 & 12 & 48 & 72 & 36 & 168 \\
4 & 0 & 12 & 84 & 240 & 384 & 720 \\
5 & 0 & 12 & 144 & 600 & 2220 & 2976 \\
6 & 0 & 12 & 252 & 1344 & 10488 & 12096 \\
7 & 0 & 12 & 456 & 2856 & 45444 & 48768 \\
8 & 0 & 12 & 852 & 5904 & 189072 & 195840 \\
\bottomrule
\end{tabular}
\caption{$M(H_4^{n};x,y)$ coefficients for $n=1,\dots,8$ (only $i\le j$ pairs).}
\label{tab:MH4n-coeffs}
\end{table}

\section*{Numerical validation}\label{sec:numerical}

In this section we validate the explicit formulation of the M-polynomial in Theorem~\ref{thm:Mpoly} and the operator–based extraction of degree–based indices summarized in Table~\ref{tab:indices}. For each family $H_p^n$ with $p\in\{3,4,5\}$ and $n=1,\dots,8$, we proceed as follows:
\begin{enumerate}
  \item We \emph{enumerate all regular states} (assignments of $n$ discs to $p$ pegs), compute the \emph{top disc} on each peg, and generate \emph{all legal moves} (top disc to an empty peg or onto a larger top disc).
  \item From the generated adjacency, we compute vertex degrees and the set of unordered edges. We then tabulate the \emph{edge multiplicities by degree pairs} $m_{i,j}(H_p^n)$ for $i\le j$.
  \item We compute the \( M \)-polynomial and evaluate the corresponding invariants 
\( M(1,1) \), \( M_1 \), \( M_2 \), \( MM_2 \), \( SSD \), \( H \), \( ISI \), \( A \), and \( F \) 
as reported in Table~\ref{tab:indices}. 
The value \( M(1,1) \) is compared with the exact number of edges \( E(H_p^n) \) 
to ensure consistency. 
We then compare all evaluated invariants with the exact results obtained using 
an algorithm that performs a complete enumeration of the graph 
and extracts the desired quantities.

\end{enumerate}

For each family we list $|E|$, $M_1$, $M_2$, $MM_2$, $SSD$, $H$, $ISI$, $A$, and $F$, computed directly from the edge degree–pair counts $m_{i,j}$:
\begin{itemize}
  \item $H_3^n$ has degree spectrum $\{2,3\}$; all indices reduce to combinations of $m_{2,2},m_{2,3},m_{3,3}$.
  \item $H_4^n$ has degree spectrum $\{3,5,6\}$; the nonzero $m_{i,j}$ appear only for $$(i,j)\in\{(3,3),(3,5),(5,5),(5,6),(6,6)\}.$$
  \item $H_5^n$ (with degree spectrum $\{4,7,9,10\}$) is included to illustrate the generality of the pipeline and to stress the scalability of the operator extraction once $m_{i,j}$ are known.
\end{itemize}

Tables~\ref{tab:indices-H3n}, \ref{tab:indices-H4n}, and~\ref{tab:indices-H5n} present the computational results.

\begin{sidewaystable}[htbp]
\centering
\small

\begin{tabular}{r|r r r r r r r r r}
\toprule
$n$ & $|E|$ & $M_1$ & $M_2$ & $MM_2$ & $SSD$ & $H$ & $ISI$ & $A$ & $F$ \\
\midrule
1 & 3 & 12 & 12 & 0.75 & 6.00 & 1.50 & 3.00 & 24.00 & 24 \\
2 & 12 & 66 & 90 & 1.67 & 25.00 & 4.40 & 16.20 & 116.34 & 186 \\
3 & 39 & 228 & 333 & 4.67 & 79.00 & 13.40 & 56.70 & 423.89 & 672 \\
4 & 120 & 714 & 1062 & 13.67 & 241.00 & 40.40 & 178.20 & 1346.53 & 2130 \\
5 & 363 & 2172 & 3249 & 40.67 & 727.00 & 121.40 & 542.70 & 4114.45 & 6504 \\
6 & 1092 & 6546 & 9810 & 121.67 & 2185.00 & 364.40 & 1636.20 & 12418.22 & 19626 \\
7 & 3279 & 19668 & 29493 & 364.67 & 6559.00 & 1093.40 & 4916.70 & 37329.52 & 58992 \\
8 & 9840 & 59034 & 88542 & 1093.67 & 19681.00 & 3280.40 & 14758.20 & 112063.41 & 177090 \\
\bottomrule
\end{tabular}
\caption{Degree-based indices for $H_{3}^{n}$, $n = 1, \dots, 8$.}
\label{tab:indices-H3n}
\vspace{.5cm}

\begin{tabular}{r|r r r r r r r r r}
\toprule
$n$ & $|E|$ & $M_1$ & $M_2$ & $MM_2$ & $SSD$ & $H$ & $ISI$ & $A$ & $F$ \\
\midrule
1 & 6 & 36 & 54 & 0.67 & 12.00 & 2.00 & 9.00 & 68.34 & 108 \\
2 & 36 & 336 & 780 & 1.76 & 75.20 & 7.80 & 82.50 & 919.92 & 1608 \\
3 & 168 & 1800 & 4836 & 6.12 & 341.60 & 31.69 & 446.86 & 5998.63 & 9792 \\
4 & 720 & 8184 & 23304 & 22.83 & 1451.20 & 127.44 & 2039.05 & 29555.77 & 46896 \\
5 & 2976 & 34776 & 101700 & 88.23 & 5975.20 & 510.89 & 8678.86 & 130380.57 & 204048 \\
6 & 12096 & 143256 & 424368 & 347.01 & 24240.00 & 2045.76 & 35781.95 & 546983.84 & 850128 \\
7 & 48768 & 581400 & 1733244 & 1376.57 & 97634.40 & 8187.47 & 145283.59 & 2240116.56 & 3469392 \\
8 & 195840 & 2342424 & 7005192 & 5483.68 & 391880.00 & 32758.85 & 585470.32 & 9066198.38 & 14016336 \\
\bottomrule
\end{tabular}
\caption{Degree-based indices for $H_{4}^{n}$, $n = 1, \dots, 8$.}
\label{tab:indices-H4n}

\vspace{.5cm}

\begin{tabular}{r|r r r r r r r r r}
\toprule
$n$ & $|E|$ & $M_1$ & $M_2$ & $MM_2$ & $SSD$ & $H$ & $ISI$ & $A$ & $F$ \\
\midrule
1 & 10 & 80 & 160 & 0.63 & 20.00 & 2.50 & 20.00 & 189.63 & 320 \\
2 & 80 & 1060 & 3500 & 1.94 & 166.43 & 12.21 & 260.91 & 4687.28 & 7180 \\
3 & 490 & 7880 & 31870 & 8.04 & 997.86 & 61.85 & 1954.66 & 47848.32 & 64640 \\
4 & 2720 & 48100 & 213740 & 36.16 & 5492.52 & 310.81 & 11973.94 & 340076.03 & 430780 \\
5 & 14410 & 268280 & 1252870 & 170.22 & 28975.00 & 1557.83 & 66909.00 & 2056518.36 & 2516720 \\
6 & 74480 & 1427860 & 6857180 & 821.77 & 149419.76 & 7799.23 & 356466.70 & 11456712.90 & 13749580 \\
7 & 379690 & 7404680 & 36143590 & 4025.16 & 760754.43 & 39023.97 & 1849636.50 & 61017378.45 & 72398240 \\
8 & 1920320 & 37831300 & 186448940 & 19882.77 & 3844775.57 & 195198.92 & 9453118.41 & 316716078.91 & 373244380 \\
\bottomrule
\end{tabular}
\caption{Degree-based indices for $H_{5}^{n}$, $n = 1, \dots, 8$.}
\label{tab:indices-H5n}
\end{sidewaystable}
We also implement computations for fixed values of \(n\) while varying \(p\). 
Thus, to complement our vertical exploration (fixed \(p\) and increasing \(n\)), 
we fix \(n \in \{1,2,3\}\) and vary \(p = 1,\dots,8\). 
Tables~\ref{tab:indices-n1-p1to8}, \ref{tab:indices-n2-p1to8}, and \ref{tab:indices-n3-p1to8} 
report the nine degree-based indices for \(H_p^n\) in this regime and confirm 
consistency with Theorem~\ref{thm:Mpoly} and Table~\ref{tab:indices}.

\begin{remark}[Analytic notes on thresholds and degree spectra (fixed $n$, varying $p$).]
    Let $r=\min\{n,p\}$ and recall that the vertex degree depends only on occupancy $\mu$ via
\[
\deg(s)=f_{p,n}(\mu)=\binom{p}{2}-\binom{p-\mu}{2},\qquad \mu\in[r].
\]
Hence the set of realized degrees is
\[
\Gamma_p^n=\big\{\,f_{p,n}(1),f_{p,n}(2),\dots,f_{p,n}(r)\,\big\},
\]
with the \emph{top-degree collision} when $r=p$:
$f_{p,n}(p-1)=f_{p,n}(p)=\Delta(H_p^n)$.

\begin{itemize}
\item Threshold $p=1$: No legal move is ever possible (only one peg), hence $E(H_1^n)=0$ for all $n$ and all degrees are $0$.

\item Threshold $p=2$: There are always moves (top disc can shuttle between the two pegs), and
\[
\Gamma_2^n=\{\,1\,\}\quad(\text{since }r=\min\{n,2\}\ge 1,\ f_{2,n}(1)=1,\ f_{2,n}(2)=1).
\]
Thus $H_2^n$ is $1$–regular (a union of cycles/paths depending on $n$), and every edge contributes the same degree pair $(1,1)$.
Consequently, all degree–based indices reduce to simple multiples of $|E|$.

\item First nontrivial heterogeneity $p=3$: Degrees split into $\Gamma_3^n=\{2,3\}$ for all $n\ge 2$ (and $\{2\}$ if $n=1$).
Only the pairs $(2,2)$, $(2,3)$, $(3,3)$ may occur. This is precisely what our $H_3^n$ tables show.

\item Rich spectra for $p\ge 4$: The realized degrees are
\[
\Gamma_p^n=\{\,p-1,\,2p-3,\,3p-6,\,\dots,\,\Delta(H_p^n)\,\},\qquad
\Delta(H_p^n)=\binom{p}{2}-1 \ \ (r=p),\ \ \text{or}\ \ \binom{p}{2}-\binom{p-r}{2}.
\]
For example:
\begin{align*}
p=4:\quad & \Gamma_4^n=\{\,3,5,6\,\},\\
p=5:\quad & \Gamma_5^n=\{\,4,7,9,10\,\},\\
p=6:\quad & \Gamma_6^n=\{\,5,9,12,14,15\,\},
\end{align*}
with the last two values equal when $r=p$ (collision at $\mu=p-1,p$).
This matches the growth and diversity visible in the $H_4^n$ and $H_5^n$ tables.

\end{itemize}
\end{remark}

\begin{table}[H]
\centering

\begin{tabular}{r|r r r r r r r r r}
\toprule
$p$ & $|E|$ & $M_1$ & $M_2$ & $MM_2$ & $SSD$ & $H$ & $ISI$ & $A$ & $F$ \\
\midrule
1 & 0 & 0 & 0 & 0.00 & 0.00 & 0.00 & 0.00 & 0.00 & 0 \\
2 & 1 & 2 & 1 & 1.00 & 2.00 & 1.00 & 0.50 & 0.00 & 2 \\
3 & 3 & 12 & 12 & 0.75 & 6.00 & 1.50 & 3.00 & 24.00 & 24 \\
4 & 6 & 36 & 54 & 0.67 & 12.00 & 2.00 & 9.00 & 68.34 & 108 \\
5 & 10 & 80 & 160 & 0.63 & 20.00 & 2.50 & 20.00 & 189.63 & 320 \\
6 & 15 & 150 & 375 & 0.60 & 30.00 & 3.00 & 37.50 & 457.76 & 750 \\
7 & 21 & 252 & 756 & 0.58 & 42.00 & 3.50 & 63.00 & 979.78 & 1512 \\
8 & 28 & 392 & 1372 & 0.57 & 56.00 & 4.00 & 98.00 & 1906.35 & 2744 \\
\bottomrule
\end{tabular}
\caption{Degree-based indices for $H_{p}^{1}$, $p = 1, \dots, 8$.}
\label{tab:indices-n1-p1to8}
\end{table}
\begin{table}[H]
\centering

\begin{tabular}{r|r r r r r r r r r}
\toprule
$p$ & $|E|$ & $M_1$ & $M_2$ & $MM_2$ & $SSD$ & $H$ & $ISI$ & $A$ & $F$ \\
\midrule
1 & 0 & 0 & 0 & 0.00 & 0.00 & 0.00 & 0.00 & 0.00 & 0 \\
2 & 2 & 4 & 2 & 2.00 & 4.00 & 2.00 & 1.00 & 0.00 & 4 \\
3 & 12 & 66 & 90 & 1.67 & 25.00 & 4.40 & 16.20 & 116.34 & 186 \\
4 & 36 & 336 & 780 & 1.76 & 75.20 & 7.80 & 82.50 & 919.92 & 1608 \\
5 & 80 & 1060 & 3500 & 1.94 & 166.43 & 12.21 & 260.91 & 4687.28 & 7180 \\
6 & 150 & 2580 & 11070 & 2.15 & 310.67 & 17.62 & 636.43 & 17151.59 & 22620 \\
7 & 252 & 5334 & 28182 & 2.37 & 519.91 & 24.03 & 1318.06 & 50081.20 & 57414 \\
8 & 392 & 9856 & 61880 & 2.60 & 806.15 & 31.45 & 2438.80 & 124554.21 & 125776 \\
\bottomrule
\end{tabular}
\caption{Degree-based indices for $H_{p}^{2}$, $p = 1, \dots, 8$.}
\label{tab:indices-n2-p1to8}
\end{table}
\begin{table}[H]
\centering
\begin{tabular}{r|r r r r r r r r r}
\toprule
$p$ & $|E|$ & $M_1$ & $M_2$ & $MM_2$ & $SSD$ & $H$ & $ISI$ & $A$ & $F$ \\
\midrule
1 & 0 & 0 & 0 & 0.00 & 0.00 & 0.00 & 0.00 & 0.00 & 0 \\
2 & 4 & 8 & 4 & 4.00 & 8.00 & 4.00 & 2.00 & 0.00 & 8 \\
3 & 39 & 228 & 333 & 4.67 & 79.00 & 13.40 & 56.70 & 423.89 & 672 \\
4 & 168 & 1800 & 4836 & 6.12 & 341.60 & 31.69 & 446.86 & 5998.63 & 9792 \\
5 & 490 & 7880 & 31870 & 8.04 & 997.86 & 61.85 & 1954.66 & 47848.32 & 64640 \\
6 & 1140 & 24720 & 135000 & 10.34 & 2320.67 & 106.90 & 6132.86 & 246376.18 & 273720 \\
7 & 2289 & 62748 & 433419 & 13.00 & 4655.00 & 169.86 & 15574.64 & 941151.50 & 877968 \\
8 & 4144 & 137648 & 1152088 & 16.01 & 8417.85 & 253.71 & 34183.57 & 2913039.64 & 2331392 \\
\bottomrule
\end{tabular}
\caption{Degree-based indices for $H_{p}^{3}$, $p = 1, \dots, 8$.}
\label{tab:indices-n3-p1to8}
\end{table}

\subsection*{Connections to Known Integer Sequences (OEIS)}
Table \ref{tab:oeis}  summarizes several striking correspondences between
degree-based topological indices of the Hanoi graphs studied in this work
and well–known sequences from the On-Line Encyclopedia of Integer Sequences
(OEIS)~\cite{OEIS}.  
For each index and each family of graphs, we list the first eight values
computed via the M-polynomial formulas established in this paper, together
with the OEIS entry that seems to match the sequence.  
The column \textit{Ment.} (“mentioned”) indicates whether the corresponding
OEIS page explicitly refers to the same graph invariant.

All correspondences were checked for the first $30$ terms and appear to hold
consistently.  
Of course, a rigorous mathematical identification requires proving that the
closed-form expressions of the corresponding indices are \emph{exactly}
equal to the OEIS formulas.  
Establishing such identities may constitute an interesting direction for
future research.

\begin{remark}
The first value listed in OEIS entry \texttt{A277105} is incorrect:
the page reports $\texttt{A277105}(1)=9$, whereas the correct value is
$MM_{2}(H_{3}^{1})=12$, as confirmed in Tables~\ref{tab:indices-H3n}
and~\ref{tab:oeis}.  
This discrepancy can be easily verified by direct computation.
\end{remark}
 
\begin{table}[H]
\centering
\resizebox{\columnwidth}{!}{
\begin{tabular}{@{}l|rrrrrrrrlcc@{}}
\toprule
$k$                                           & 1  & 2  & 3   & 4    & 5    & 6    & 7     & 8     & OEIS                    & \multicolumn{1}{l}{Condition} & \multicolumn{1}{l}{Ment.} \\ \midrule
$M_{1}(H_{3}^{k})$                            & 12 & 66 & 228 & 714  & 2172 & 6546 & 19668 & 59034 & $\texttt{A277104}(k)$   &                               & yes                           \\
$M_{2}(H_{3}^{k})$                            & 12 & 90 & 333 & 1062 & 3249 & 9810 & 29493 & 88542 & $\texttt{A277105}(k)$   & $k\geq2$                      & yes                           \\
$\left\lfloor MM_{2}(H_{3}^{k})\right\rfloor$ & 0  & 1  & 4   & 13   & 40   & 121  & 364   & 1093  & $\texttt{A003462}(k-1)$ &                               & no                            \\
$SSD(H_{3}^{k})$                              & 6  & 25 & 79  & 241  & 727  & 2185 & 6559  & 19681 & $\texttt{A058481}(k+1)$ & $k\geq2$                      & no                            \\
$\left\lfloor H(H_{3}^{k})\right\rfloor$      & 1  & 4  & 13  & 40   & 121  & 364  & 1093  & 3280  & $\texttt{A003462}(k)$   &                               & no                            \\
$M_{1}(H_{k}^{1})$                            & 0  & 2  & 12  & 36   & 80   & 150  & 252   & 392   & $\texttt{A011379}(k-1)$ &                               & yes                           \\
$M_{2}(H_{k}^{1})$                            & 0  & 1  & 12  & 54   & 160  & 375  & 756   & 1372  & $\texttt{A019582}(k)$   &                               & yes                           \\
$SSD(H_{k}^{1})$                              & 0  & 2  & 6   & 12   & 20   & 30   & 42    & 56    & $\texttt{A002378}(k-1)$ &                               & no                            \\
$F(H_{k}^{1})$                                & 0  & 2  & 24  & 108  & 320  & 750  & 1512  & 2744  & $\texttt{A179824}(k)$   & $k\geq2$                      & no                            \\
$\left\lfloor MM_{2}(H_{k}^{3})\right\rfloor$ & 0  & 4  & 4   & 6    & 8    & 10   & 13    & 16    & $\texttt{A008748}(k+1)$ & $k\geq3$                      & no                            \\
$\left\lceil MM_{2}(H_{k}^{3})\right\rceil$   & 0  & 4  & 5   & 7    & 9    & 11   & 14    & 17    & $\texttt{A008750}(k+2)$ & $k\geq2$                      & no                            \\ \bottomrule
\end{tabular}
}
\caption{Correspondences between degree-based topological indices of
$H_{p}^{n}$ and known integer sequences in OEIS.  
For each index, we list its first eight values, the associated OEIS entry
(if any), the range of validity (when applicable), and whether the index is 
explicitly mentioned on the OEIS page.}
\label{tab:oeis}
\end{table}

\section{Conclusion}\label{sec:conclusion}

In this work, we derived the first complete and explicit formulation of the M-polynomial of the
generalized Hanoi graphs $H_p^n$, a natural yet previously unexplored family of graphs arising
from the Tower of Hanoi problem with an arbitrary number of pegs $p $ and discs $n$.
Our contribution intersects combinatorial enumeration, algebraic graph theory, and the theory of
degree-based topological indices, significantly expanding the applicability of the M-polynomial
framework introduced by Deutsch and Klavžar \cite{EmericMploy}.

The central achievement of this paper is a detailed combinatorial decomposition of the edge set of
$H_p^n$ based on occupancy transitions, distinguishing moves between singleton, multiton, empty,
and occupied pegs. This structural analysis yields fully explicit formulas for all coefficients of the
M-polynomial. Remarkably, the entire polynomial can be expressed using only falling factorials and
(associated) Stirling numbers of the second kind, ensuring  computational
tractability for arbitrary $p$ and $n$.

Our approach simultaneously characterizes the degree structure of the generalized Hanoi graphs,
the adjacency relations between degree classes, and the combinatorial constraints imposed by the
Tower of Hanoi rules. These insights culminate in Theorem~\ref{thm:Mpoly}, which provides the
first exact M-polynomial for $H_p^n$ in full generality. Furthermore, substituting this expression
into the differential identities of Table~\ref{tab:indices} yields closed formulas for all classical
degree-based indices (Zagreb, Randić, harmonic, symmetric division, forgotten index, etc.).
The correctness of these results has been validated through explicit enumeration for representative
instances, confirming the structural predictions and functional extraction properties of our formulas.

\bibliographystyle{plain}
\bibliography{bibou}

\end{document}

%% file: small_graphs.tex
\begin{figure}[H]
    \centering
    
 \resizebox{\columnwidth}{!}{
\tikzset{every picture/.style={line width=0.75pt}} 

\begin{tikzpicture}[x=0.75pt,y=0.75pt,yscale=-1,xscale=1]

\draw    (324.5,55.64) -- (338.64,69.78) ;
\draw    (352.78,55.64) -- (338.64,69.78) ;
\draw    (338.64,41.5) -- (324.5,55.64) ;
\draw    (338.64,41.5) -- (352.78,55.64) ;
\draw  [fill={rgb, 255:red, 0; green, 0; blue, 0 }  ,fill opacity=1 ] (351.02,53.87) .. controls (351.99,52.9) and (353.58,52.9) .. (354.55,53.87) .. controls (355.53,54.85) and (355.53,56.43) .. (354.55,57.41) .. controls (353.58,58.39) and (351.99,58.39) .. (351.02,57.41) .. controls (350.04,56.43) and (350.04,54.85) .. (351.02,53.87) -- cycle ;
\draw  [fill={rgb, 255:red, 0; green, 0; blue, 0 }  ,fill opacity=1 ] (336.87,39.73) .. controls (337.85,38.76) and (339.43,38.76) .. (340.41,39.73) .. controls (341.39,40.71) and (341.39,42.29) .. (340.41,43.27) .. controls (339.43,44.24) and (337.85,44.24) .. (336.87,43.27) .. controls (335.9,42.29) and (335.9,40.71) .. (336.87,39.73) -- cycle ;
\draw  [fill={rgb, 255:red, 0; green, 0; blue, 0 }  ,fill opacity=1 ] (336.87,68.02) .. controls (337.85,67.04) and (339.43,67.04) .. (340.41,68.02) .. controls (341.39,68.99) and (341.39,70.58) .. (340.41,71.55) .. controls (339.43,72.53) and (337.85,72.53) .. (336.87,71.55) .. controls (335.9,70.58) and (335.9,68.99) .. (336.87,68.02) -- cycle ;
\draw  [fill={rgb, 255:red, 0; green, 0; blue, 0 }  ,fill opacity=1 ] (322.73,53.87) .. controls (323.71,52.9) and (325.29,52.9) .. (326.27,53.87) .. controls (327.24,54.85) and (327.24,56.43) .. (326.27,57.41) .. controls (325.29,58.39) and (323.71,58.39) .. (322.73,57.41) .. controls (321.76,56.43) and (321.76,54.85) .. (322.73,53.87) -- cycle ;
\draw    (324.5,55.64) -- (352.78,55.64) ;
\draw    (338.64,41.5) -- (338.64,69.78) ;
\draw    (374.5,105.64) -- (388.64,119.78) ;
\draw    (402.78,105.64) -- (388.64,119.78) ;
\draw    (388.64,91.5) -- (374.5,105.64) ;
\draw    (388.64,91.5) -- (402.78,105.64) ;
\draw  [fill={rgb, 255:red, 0; green, 0; blue, 0 }  ,fill opacity=1 ] (401.02,103.87) .. controls (401.99,102.9) and (403.58,102.9) .. (404.55,103.87) .. controls (405.53,104.85) and (405.53,106.43) .. (404.55,107.41) .. controls (403.58,108.39) and (401.99,108.39) .. (401.02,107.41) .. controls (400.04,106.43) and (400.04,104.85) .. (401.02,103.87) -- cycle ;
\draw  [fill={rgb, 255:red, 0; green, 0; blue, 0 }  ,fill opacity=1 ] (386.87,89.73) .. controls (387.85,88.76) and (389.43,88.76) .. (390.41,89.73) .. controls (391.39,90.71) and (391.39,92.29) .. (390.41,93.27) .. controls (389.43,94.24) and (387.85,94.24) .. (386.87,93.27) .. controls (385.9,92.29) and (385.9,90.71) .. (386.87,89.73) -- cycle ;
\draw  [fill={rgb, 255:red, 0; green, 0; blue, 0 }  ,fill opacity=1 ] (386.87,118.02) .. controls (387.85,117.04) and (389.43,117.04) .. (390.41,118.02) .. controls (391.39,118.99) and (391.39,120.58) .. (390.41,121.55) .. controls (389.43,122.53) and (387.85,122.53) .. (386.87,121.55) .. controls (385.9,120.58) and (385.9,118.99) .. (386.87,118.02) -- cycle ;
\draw  [fill={rgb, 255:red, 0; green, 0; blue, 0 }  ,fill opacity=1 ] (372.73,103.87) .. controls (373.71,102.9) and (375.29,102.9) .. (376.27,103.87) .. controls (377.24,104.85) and (377.24,106.43) .. (376.27,107.41) .. controls (375.29,108.39) and (373.71,108.39) .. (372.73,107.41) .. controls (371.76,106.43) and (371.76,104.85) .. (372.73,103.87) -- cycle ;
\draw    (374.5,105.64) -- (402.78,105.64) ;
\draw    (388.64,91.5) -- (388.64,119.78) ;
\draw    (274.5,105.64) -- (288.64,119.78) ;
\draw    (302.78,105.64) -- (288.64,119.78) ;
\draw    (288.64,91.5) -- (274.5,105.64) ;
\draw    (288.64,91.5) -- (302.78,105.64) ;
\draw  [fill={rgb, 255:red, 0; green, 0; blue, 0 }  ,fill opacity=1 ] (301.02,103.87) .. controls (301.99,102.9) and (303.58,102.9) .. (304.55,103.87) .. controls (305.53,104.85) and (305.53,106.43) .. (304.55,107.41) .. controls (303.58,108.39) and (301.99,108.39) .. (301.02,107.41) .. controls (300.04,106.43) and (300.04,104.85) .. (301.02,103.87) -- cycle ;
\draw  [fill={rgb, 255:red, 0; green, 0; blue, 0 }  ,fill opacity=1 ] (286.87,89.73) .. controls (287.85,88.76) and (289.43,88.76) .. (290.41,89.73) .. controls (291.39,90.71) and (291.39,92.29) .. (290.41,93.27) .. controls (289.43,94.24) and (287.85,94.24) .. (286.87,93.27) .. controls (285.9,92.29) and (285.9,90.71) .. (286.87,89.73) -- cycle ;
\draw  [fill={rgb, 255:red, 0; green, 0; blue, 0 }  ,fill opacity=1 ] (286.87,118.02) .. controls (287.85,117.04) and (289.43,117.04) .. (290.41,118.02) .. controls (291.39,118.99) and (291.39,120.58) .. (290.41,121.55) .. controls (289.43,122.53) and (287.85,122.53) .. (286.87,121.55) .. controls (285.9,120.58) and (285.9,118.99) .. (286.87,118.02) -- cycle ;
\draw  [fill={rgb, 255:red, 0; green, 0; blue, 0 }  ,fill opacity=1 ] (272.73,103.87) .. controls (273.71,102.9) and (275.29,102.9) .. (276.27,103.87) .. controls (277.24,104.85) and (277.24,106.43) .. (276.27,107.41) .. controls (275.29,108.39) and (273.71,108.39) .. (272.73,107.41) .. controls (271.76,106.43) and (271.76,104.85) .. (272.73,103.87) -- cycle ;
\draw    (274.5,105.64) -- (302.78,105.64) ;
\draw    (288.64,91.5) -- (288.64,119.78) ;
\draw    (324.5,155.64) -- (338.64,169.78) ;
\draw    (352.78,155.64) -- (338.64,169.78) ;
\draw    (338.64,141.5) -- (324.5,155.64) ;
\draw    (338.64,141.5) -- (352.78,155.64) ;
\draw  [fill={rgb, 255:red, 0; green, 0; blue, 0 }  ,fill opacity=1 ] (351.02,153.87) .. controls (351.99,152.9) and (353.58,152.9) .. (354.55,153.87) .. controls (355.53,154.85) and (355.53,156.43) .. (354.55,157.41) .. controls (353.58,158.39) and (351.99,158.39) .. (351.02,157.41) .. controls (350.04,156.43) and (350.04,154.85) .. (351.02,153.87) -- cycle ;
\draw  [fill={rgb, 255:red, 0; green, 0; blue, 0 }  ,fill opacity=1 ] (336.87,139.73) .. controls (337.85,138.76) and (339.43,138.76) .. (340.41,139.73) .. controls (341.39,140.71) and (341.39,142.29) .. (340.41,143.27) .. controls (339.43,144.24) and (337.85,144.24) .. (336.87,143.27) .. controls (335.9,142.29) and (335.9,140.71) .. (336.87,139.73) -- cycle ;
\draw  [fill={rgb, 255:red, 0; green, 0; blue, 0 }  ,fill opacity=1 ] (336.87,168.02) .. controls (337.85,167.04) and (339.43,167.04) .. (340.41,168.02) .. controls (341.39,168.99) and (341.39,170.58) .. (340.41,171.55) .. controls (339.43,172.53) and (337.85,172.53) .. (336.87,171.55) .. controls (335.9,170.58) and (335.9,168.99) .. (336.87,168.02) -- cycle ;
\draw  [fill={rgb, 255:red, 0; green, 0; blue, 0 }  ,fill opacity=1 ] (322.73,153.87) .. controls (323.71,152.9) and (325.29,152.9) .. (326.27,153.87) .. controls (327.24,154.85) and (327.24,156.43) .. (326.27,157.41) .. controls (325.29,158.39) and (323.71,158.39) .. (322.73,157.41) .. controls (321.76,156.43) and (321.76,154.85) .. (322.73,153.87) -- cycle ;
\draw    (324.5,155.64) -- (352.78,155.64) ;
\draw    (338.64,141.5) -- (338.64,169.78) ;
\draw    (352.78,55.64) -- (352.78,155.64) ;
\draw    (324.5,55.64) -- (324.5,155.64) ;
\draw    (352.78,55.64) -- (374.5,105.64) ;
\draw    (302.78,105.64) -- (324.5,155.64) ;
\draw    (324.5,55.64) -- (302.78,105.64) ;
\draw    (374.5,105.64) -- (352.78,155.64) ;
\draw    (388.64,91.5) -- (288.64,91.5) ;
\draw    (388.64,119.78) -- (288.64,119.78) ;
\draw    (388.64,91.5) -- (338.64,69.78) ;
\draw    (338.64,141.5) -- (288.64,119.78) ;
\draw    (388.64,119.78) -- (338.64,141.5) ;
\draw    (338.64,69.78) -- (288.64,91.5) ;

\draw  [fill={rgb, 255:red, 0; green, 0; blue, 0 }  ,fill opacity=1 ] (114.83,50.14) .. controls (114.83,48.76) and (115.95,47.64) .. (117.33,47.64) .. controls (118.71,47.64) and (119.83,48.76) .. (119.83,50.14) .. controls (119.83,51.52) and (118.71,52.64) .. (117.33,52.64) .. controls (115.95,52.64) and (114.83,51.52) .. (114.83,50.14) -- cycle ;
\draw    (107.33,67.46) -- (127.33,67.46) ;
\draw    (117.33,50.14) -- (127.33,67.46) ;
\draw    (117.33,50.14) -- (107.33,67.46) ;
\draw  [fill={rgb, 255:red, 0; green, 0; blue, 0 }  ,fill opacity=1 ] (124.83,67.46) .. controls (124.83,66.08) and (125.95,64.96) .. (127.33,64.96) .. controls (128.71,64.96) and (129.83,66.08) .. (129.83,67.46) .. controls (129.83,68.84) and (128.71,69.96) .. (127.33,69.96) .. controls (125.95,69.96) and (124.83,68.84) .. (124.83,67.46) -- cycle ;
\draw  [fill={rgb, 255:red, 0; green, 0; blue, 0 }  ,fill opacity=1 ] (104.83,67.46) .. controls (104.83,66.08) and (105.95,64.96) .. (107.33,64.96) .. controls (108.71,64.96) and (109.83,66.08) .. (109.83,67.46) .. controls (109.83,68.84) and (108.71,69.96) .. (107.33,69.96) .. controls (105.95,69.96) and (104.83,68.84) .. (104.83,67.46) -- cycle ;
\draw  [fill={rgb, 255:red, 0; green, 0; blue, 0 }  ,fill opacity=1 ] (134.83,84.78) .. controls (134.83,83.4) and (135.95,82.28) .. (137.33,82.28) .. controls (138.71,82.28) and (139.83,83.4) .. (139.83,84.78) .. controls (139.83,86.16) and (138.71,87.28) .. (137.33,87.28) .. controls (135.95,87.28) and (134.83,86.16) .. (134.83,84.78) -- cycle ;
\draw    (127.33,102.1) -- (147.33,102.1) ;
\draw    (137.33,84.78) -- (147.33,102.1) ;
\draw    (137.33,84.78) -- (127.33,102.1) ;
\draw  [fill={rgb, 255:red, 0; green, 0; blue, 0 }  ,fill opacity=1 ] (144.83,102.1) .. controls (144.83,100.72) and (145.95,99.6) .. (147.33,99.6) .. controls (148.71,99.6) and (149.83,100.72) .. (149.83,102.1) .. controls (149.83,103.48) and (148.71,104.6) .. (147.33,104.6) .. controls (145.95,104.6) and (144.83,103.48) .. (144.83,102.1) -- cycle ;
\draw  [fill={rgb, 255:red, 0; green, 0; blue, 0 }  ,fill opacity=1 ] (124.83,102.1) .. controls (124.83,100.72) and (125.95,99.6) .. (127.33,99.6) .. controls (128.71,99.6) and (129.83,100.72) .. (129.83,102.1) .. controls (129.83,103.48) and (128.71,104.6) .. (127.33,104.6) .. controls (125.95,104.6) and (124.83,103.48) .. (124.83,102.1) -- cycle ;
\draw  [fill={rgb, 255:red, 0; green, 0; blue, 0 }  ,fill opacity=1 ] (94.83,84.78) .. controls (94.83,83.4) and (95.95,82.28) .. (97.33,82.28) .. controls (98.71,82.28) and (99.83,83.4) .. (99.83,84.78) .. controls (99.83,86.16) and (98.71,87.28) .. (97.33,87.28) .. controls (95.95,87.28) and (94.83,86.16) .. (94.83,84.78) -- cycle ;
\draw    (87.33,102.1) -- (107.33,102.1) ;
\draw    (97.33,84.78) -- (107.33,102.1) ;
\draw    (97.33,84.78) -- (87.33,102.1) ;
\draw  [fill={rgb, 255:red, 0; green, 0; blue, 0 }  ,fill opacity=1 ] (104.83,102.1) .. controls (104.83,100.72) and (105.95,99.6) .. (107.33,99.6) .. controls (108.71,99.6) and (109.83,100.72) .. (109.83,102.1) .. controls (109.83,103.48) and (108.71,104.6) .. (107.33,104.6) .. controls (105.95,104.6) and (104.83,103.48) .. (104.83,102.1) -- cycle ;
\draw  [fill={rgb, 255:red, 0; green, 0; blue, 0 }  ,fill opacity=1 ] (84.83,102.1) .. controls (84.83,100.72) and (85.95,99.6) .. (87.33,99.6) .. controls (88.71,99.6) and (89.83,100.72) .. (89.83,102.1) .. controls (89.83,103.48) and (88.71,104.6) .. (87.33,104.6) .. controls (85.95,104.6) and (84.83,103.48) .. (84.83,102.1) -- cycle ;
\draw    (127.33,67.46) -- (137.33,84.78) ;
\draw    (107.33,102.1) -- (127.33,102.1) ;
\draw    (107.33,67.46) -- (97.33,84.78) ;
\draw  [fill={rgb, 255:red, 0; green, 0; blue, 0 }  ,fill opacity=1 ] (154.83,119.64) .. controls (154.83,118.26) and (155.95,117.14) .. (157.33,117.14) .. controls (158.71,117.14) and (159.83,118.26) .. (159.83,119.64) .. controls (159.83,121.02) and (158.71,122.14) .. (157.33,122.14) .. controls (155.95,122.14) and (154.83,121.02) .. (154.83,119.64) -- cycle ;
\draw    (147.33,136.96) -- (167.33,136.96) ;
\draw    (157.33,119.64) -- (167.33,136.96) ;
\draw    (157.33,119.64) -- (147.33,136.96) ;
\draw  [fill={rgb, 255:red, 0; green, 0; blue, 0 }  ,fill opacity=1 ] (164.83,136.96) .. controls (164.83,135.58) and (165.95,134.46) .. (167.33,134.46) .. controls (168.71,134.46) and (169.83,135.58) .. (169.83,136.96) .. controls (169.83,138.34) and (168.71,139.46) .. (167.33,139.46) .. controls (165.95,139.46) and (164.83,138.34) .. (164.83,136.96) -- cycle ;
\draw  [fill={rgb, 255:red, 0; green, 0; blue, 0 }  ,fill opacity=1 ] (144.83,136.96) .. controls (144.83,135.58) and (145.95,134.46) .. (147.33,134.46) .. controls (148.71,134.46) and (149.83,135.58) .. (149.83,136.96) .. controls (149.83,138.34) and (148.71,139.46) .. (147.33,139.46) .. controls (145.95,139.46) and (144.83,138.34) .. (144.83,136.96) -- cycle ;
\draw  [fill={rgb, 255:red, 0; green, 0; blue, 0 }  ,fill opacity=1 ] (174.83,154.28) .. controls (174.83,152.9) and (175.95,151.78) .. (177.33,151.78) .. controls (178.71,151.78) and (179.83,152.9) .. (179.83,154.28) .. controls (179.83,155.66) and (178.71,156.78) .. (177.33,156.78) .. controls (175.95,156.78) and (174.83,155.66) .. (174.83,154.28) -- cycle ;
\draw    (167.33,171.6) -- (187.33,171.6) ;
\draw    (177.33,154.28) -- (187.33,171.6) ;
\draw    (177.33,154.28) -- (167.33,171.6) ;
\draw  [fill={rgb, 255:red, 0; green, 0; blue, 0 }  ,fill opacity=1 ] (184.83,171.6) .. controls (184.83,170.22) and (185.95,169.1) .. (187.33,169.1) .. controls (188.71,169.1) and (189.83,170.22) .. (189.83,171.6) .. controls (189.83,172.98) and (188.71,174.1) .. (187.33,174.1) .. controls (185.95,174.1) and (184.83,172.98) .. (184.83,171.6) -- cycle ;
\draw  [fill={rgb, 255:red, 0; green, 0; blue, 0 }  ,fill opacity=1 ] (164.83,171.6) .. controls (164.83,170.22) and (165.95,169.1) .. (167.33,169.1) .. controls (168.71,169.1) and (169.83,170.22) .. (169.83,171.6) .. controls (169.83,172.98) and (168.71,174.1) .. (167.33,174.1) .. controls (165.95,174.1) and (164.83,172.98) .. (164.83,171.6) -- cycle ;
\draw  [fill={rgb, 255:red, 0; green, 0; blue, 0 }  ,fill opacity=1 ] (134.83,154.28) .. controls (134.83,152.9) and (135.95,151.78) .. (137.33,151.78) .. controls (138.71,151.78) and (139.83,152.9) .. (139.83,154.28) .. controls (139.83,155.66) and (138.71,156.78) .. (137.33,156.78) .. controls (135.95,156.78) and (134.83,155.66) .. (134.83,154.28) -- cycle ;
\draw    (127.33,171.6) -- (147.33,171.6) ;
\draw    (137.33,154.28) -- (147.33,171.6) ;
\draw    (137.33,154.28) -- (127.33,171.6) ;
\draw  [fill={rgb, 255:red, 0; green, 0; blue, 0 }  ,fill opacity=1 ] (144.83,171.6) .. controls (144.83,170.22) and (145.95,169.1) .. (147.33,169.1) .. controls (148.71,169.1) and (149.83,170.22) .. (149.83,171.6) .. controls (149.83,172.98) and (148.71,174.1) .. (147.33,174.1) .. controls (145.95,174.1) and (144.83,172.98) .. (144.83,171.6) -- cycle ;
\draw  [fill={rgb, 255:red, 0; green, 0; blue, 0 }  ,fill opacity=1 ] (124.83,171.6) .. controls (124.83,170.22) and (125.95,169.1) .. (127.33,169.1) .. controls (128.71,169.1) and (129.83,170.22) .. (129.83,171.6) .. controls (129.83,172.98) and (128.71,174.1) .. (127.33,174.1) .. controls (125.95,174.1) and (124.83,172.98) .. (124.83,171.6) -- cycle ;
\draw    (167.33,136.96) -- (177.33,154.28) ;
\draw    (147.33,171.6) -- (167.33,171.6) ;
\draw    (147.33,136.96) -- (137.33,154.28) ;
\draw    (87.33,102.1) -- (77.33,119.42) ;
\draw    (147.33,102.1) -- (157.33,119.42) ;
\draw  [fill={rgb, 255:red, 0; green, 0; blue, 0 }  ,fill opacity=1 ] (74.83,119.84) .. controls (74.83,118.46) and (75.95,117.34) .. (77.33,117.34) .. controls (78.71,117.34) and (79.83,118.46) .. (79.83,119.84) .. controls (79.83,121.22) and (78.71,122.34) .. (77.33,122.34) .. controls (75.95,122.34) and (74.83,121.22) .. (74.83,119.84) -- cycle ;
\draw    (67.33,137.16) -- (87.33,137.16) ;
\draw    (77.33,119.84) -- (87.33,137.16) ;
\draw    (77.33,119.84) -- (67.33,137.16) ;
\draw  [fill={rgb, 255:red, 0; green, 0; blue, 0 }  ,fill opacity=1 ] (84.83,137.16) .. controls (84.83,135.78) and (85.95,134.66) .. (87.33,134.66) .. controls (88.71,134.66) and (89.83,135.78) .. (89.83,137.16) .. controls (89.83,138.54) and (88.71,139.66) .. (87.33,139.66) .. controls (85.95,139.66) and (84.83,138.54) .. (84.83,137.16) -- cycle ;
\draw  [fill={rgb, 255:red, 0; green, 0; blue, 0 }  ,fill opacity=1 ] (64.83,137.16) .. controls (64.83,135.78) and (65.95,134.66) .. (67.33,134.66) .. controls (68.71,134.66) and (69.83,135.78) .. (69.83,137.16) .. controls (69.83,138.54) and (68.71,139.66) .. (67.33,139.66) .. controls (65.95,139.66) and (64.83,138.54) .. (64.83,137.16) -- cycle ;
\draw  [fill={rgb, 255:red, 0; green, 0; blue, 0 }  ,fill opacity=1 ] (94.83,154.48) .. controls (94.83,153.1) and (95.95,151.98) .. (97.33,151.98) .. controls (98.71,151.98) and (99.83,153.1) .. (99.83,154.48) .. controls (99.83,155.86) and (98.71,156.98) .. (97.33,156.98) .. controls (95.95,156.98) and (94.83,155.86) .. (94.83,154.48) -- cycle ;
\draw    (87.33,171.8) -- (107.33,171.8) ;
\draw    (97.33,154.48) -- (107.33,171.8) ;
\draw    (97.33,154.48) -- (87.33,171.8) ;
\draw  [fill={rgb, 255:red, 0; green, 0; blue, 0 }  ,fill opacity=1 ] (104.83,171.8) .. controls (104.83,170.42) and (105.95,169.3) .. (107.33,169.3) .. controls (108.71,169.3) and (109.83,170.42) .. (109.83,171.8) .. controls (109.83,173.18) and (108.71,174.3) .. (107.33,174.3) .. controls (105.95,174.3) and (104.83,173.18) .. (104.83,171.8) -- cycle ;
\draw  [fill={rgb, 255:red, 0; green, 0; blue, 0 }  ,fill opacity=1 ] (84.83,171.8) .. controls (84.83,170.42) and (85.95,169.3) .. (87.33,169.3) .. controls (88.71,169.3) and (89.83,170.42) .. (89.83,171.8) .. controls (89.83,173.18) and (88.71,174.3) .. (87.33,174.3) .. controls (85.95,174.3) and (84.83,173.18) .. (84.83,171.8) -- cycle ;
\draw  [fill={rgb, 255:red, 0; green, 0; blue, 0 }  ,fill opacity=1 ] (54.83,154.48) .. controls (54.83,153.1) and (55.95,151.98) .. (57.33,151.98) .. controls (58.71,151.98) and (59.83,153.1) .. (59.83,154.48) .. controls (59.83,155.86) and (58.71,156.98) .. (57.33,156.98) .. controls (55.95,156.98) and (54.83,155.86) .. (54.83,154.48) -- cycle ;
\draw    (47.33,171.8) -- (67.33,171.8) ;
\draw    (57.33,154.48) -- (67.33,171.8) ;
\draw    (57.33,154.48) -- (47.33,171.8) ;
\draw  [fill={rgb, 255:red, 0; green, 0; blue, 0 }  ,fill opacity=1 ] (64.83,171.8) .. controls (64.83,170.42) and (65.95,169.3) .. (67.33,169.3) .. controls (68.71,169.3) and (69.83,170.42) .. (69.83,171.8) .. controls (69.83,173.18) and (68.71,174.3) .. (67.33,174.3) .. controls (65.95,174.3) and (64.83,173.18) .. (64.83,171.8) -- cycle ;
\draw  [fill={rgb, 255:red, 0; green, 0; blue, 0 }  ,fill opacity=1 ] (44.83,171.8) .. controls (44.83,170.42) and (45.95,169.3) .. (47.33,169.3) .. controls (48.71,169.3) and (49.83,170.42) .. (49.83,171.8) .. controls (49.83,173.18) and (48.71,174.3) .. (47.33,174.3) .. controls (45.95,174.3) and (44.83,173.18) .. (44.83,171.8) -- cycle ;
\draw    (87.33,137.16) -- (97.33,154.48) ;
\draw    (67.33,171.8) -- (87.33,171.8) ;
\draw    (67.33,137.16) -- (57.33,154.48) ;
\draw    (107.33,171.6) -- (127.33,171.6) ;

\draw (111,36) node [anchor=north west][inner sep=0.75pt]  [font=\tiny]  {$000$};
\draw (134.15,66.16) node [anchor=north west][inner sep=0.75pt]  [font=\tiny]  {$001$};
\draw (144.55,83.48) node [anchor=north west][inner sep=0.75pt]  [font=\tiny]  {$021$};
\draw (174.15,135.45) node [anchor=north west][inner sep=0.75pt]  [font=\tiny]  {$120$};
\draw (164.15,118.13) node [anchor=north west][inner sep=0.75pt]  [font=\tiny]  {$122$};
\draw (184.55,152.77) node [anchor=north west][inner sep=0.75pt]  [font=\tiny]  {$110$};
\draw (27.83,170.12) node [anchor=north west][inner sep=0.75pt]  [font=\tiny]  {$222$};
\draw (37.83,152.8) node [anchor=north west][inner sep=0.75pt]  [font=\tiny]  {$220$};
\draw (47.83,135.48) node [anchor=north west][inner sep=0.75pt]  [font=\tiny]  {$210$};
\draw (87.83,66.38) node [anchor=north west][inner sep=0.75pt]  [font=\tiny]  {$002$};
\draw (77.83,83.52) node [anchor=north west][inner sep=0.75pt]  [font=\tiny]  {$012$};
\draw (57.83,118.16) node [anchor=north west][inner sep=0.75pt]  [font=\tiny]  {$211$};
\draw (67.83,100.84) node [anchor=north west][inner sep=0.75pt]  [font=\tiny]  {$011$};
\draw (154.15,100.81) node [anchor=north west][inner sep=0.75pt]  [font=\tiny]  {$022$};
\draw (194.15,169.59) node [anchor=north west][inner sep=0.75pt]  [font=\tiny]  {$111$};
\draw (92.55,135.21) node [anchor=north west][inner sep=0.75pt]  [font=\tiny]  {$212$};
\draw (103.55,152.51) node [anchor=north west][inner sep=0.75pt]  [font=\tiny]  {$202$};
\draw (102.33,110) node [anchor=north west][inner sep=0.75pt]  [font=\tiny]  {$010$};
\draw (128.55,135.21) node [anchor=north west][inner sep=0.75pt]  [font=\tiny]  {$121$};
\draw (121.33,179.5) node [anchor=north west][inner sep=0.75pt]  [font=\tiny]  {$100$};
\draw (101.33,179.7) node [anchor=north west][inner sep=0.75pt]  [font=\tiny]  {$200$};
\draw (141.33,179.5) node [anchor=north west][inner sep=0.75pt]  [font=\tiny]  {$102$};
\draw (161.33,179.5) node [anchor=north west][inner sep=0.75pt]  [font=\tiny]  {$112$};
\draw (81.33,179.7) node [anchor=north west][inner sep=0.75pt]  [font=\tiny]  {$201$};
\draw (61.33,179.7) node [anchor=north west][inner sep=0.75pt]  [font=\tiny]  {$221$};
\draw (122.33,110) node [anchor=north west][inner sep=0.75pt]  [font=\tiny]  {$020$};
\draw (120,152.51) node [anchor=north west][inner sep=0.75pt]  [font=\tiny]  {$101$};
\draw (335,30) node [anchor=north west][inner sep=0.75pt]  [font=\tiny]  {$00$};
\draw (307,54.4) node [anchor=north west][inner sep=0.75pt]  [font=\tiny]  {$01$};
\draw (359,54.4) node [anchor=north west][inner sep=0.75pt]  [font=\tiny]  {$03$};
\draw (335,79.13) node [anchor=north west][inner sep=0.75pt]  [font=\tiny]  {$02$};
\draw (385,128.4) node [anchor=north west][inner sep=0.75pt]  [font=\tiny]  {$10$};
\draw (358,106.4) node [anchor=north west][inner sep=0.75pt]  [font=\tiny]  {$13$};
\draw (385,80) node [anchor=north west][inner sep=0.75pt]  [font=\tiny]  {$12$};
\draw (409.3,106.4) node [anchor=north west][inner sep=0.75pt]  [font=\tiny]  {$11$};
\draw (335,130) node [anchor=north west][inner sep=0.75pt]  [font=\tiny]  {$20$};
\draw (335,179.7) node [anchor=north west][inner sep=0.75pt]  [font=\tiny]  {$23$};
\draw (359,157.4) node [anchor=north west][inner sep=0.75pt]  [font=\tiny]  {$22$};
\draw (312,157.4) node [anchor=north west][inner sep=0.75pt]  [font=\tiny]  {$21$};
\draw (282,128.4) node [anchor=north west][inner sep=0.75pt]  [font=\tiny]  {$30$};
\draw (259.7,106.4) node [anchor=north west][inner sep=0.75pt]  [font=\tiny]  {$33$};
\draw (282,80) node [anchor=north west][inner sep=0.75pt]  [font=\tiny]  {$32$};
\draw (309.3,106.4) node [anchor=north west][inner sep=0.75pt]  [font=\tiny]  {$31$};
\draw (105.49,207.98) node [anchor=north west][inner sep=0.75pt]  []  {$H_{3}^{3}$};
\draw (329,207.98) node [anchor=north west][inner sep=0.75pt]  []  {$H_{4}^{2}$};

\end{tikzpicture}

}
     \caption{Examples of Hanoi graphs \(H_3^3\) and \(H_4^2\).  
    Each vertex represents a regular configuration of discs, and edges represent single legal moves.}
    \label{fig:small_graphs}
\end{figure}